\RequirePackage{fix-cm}
\documentclass[smallextended]{svjour3}       % onecolumn (second format)

\smartqed  % flush right qed marks, e.g. at end of proof

\usepackage[margin=1.15in]{geometry}
\usepackage{subfigure}
\usepackage{graphicx}
\usepackage{amsfonts}
\usepackage{amsmath}
\usepackage{listings,chngcntr}
\usepackage{xcolor}
\usepackage{wrapfig}

\definecolor{codegreen}{rgb}{0,0.6,0}
\definecolor{codegray}{rgb}{0.5,0.5,0.5}
\definecolor{codepurple}{rgb}{0.58,0,0.82}
\definecolor{backcolour}{rgb}{0.95,0.95,0.92}

\lstdefinestyle{mystyle}{
    backgroundcolor=\color{backcolour},
    commentstyle=\color{codegreen},
    keywordstyle=\color{magenta},
    numberstyle=\tiny\color{codegray},
    stringstyle=\color{codepurple},
    basicstyle=\ttfamily\footnotesize,
    breakatwhitespace=false,
    breaklines=true,
    captionpos=b,
    keepspaces=true,
    numbers=left,
    numbersep=5pt,
    showspaces=false,
    showstringspaces=false,
    showtabs=false,
    tabsize=2
}

\usepackage{fontawesome5}
\usepackage{hyperref}
\usepackage[section]{algorithm}
\usepackage{algorithmic}
\numberwithin{equation}{section}
\begin{document}
%\title{Using Python to teach a basic forecasting course}
\title{A basic time series forecasting course with Python %\thanks{Grants or other notes
}
%%%%%%%%%%%%%%%%
%{{\color{blue}{\textsf{Python}}}}
%{\textbf{\includegraphics[width=2cm]{Python-logo-master-v3-TM.png}}


\author{Alain Zemkoho %        \and
        %Second Author %etc.
}


\institute{Alain Zemkoho \at
             School of Mathematical Sciences \\
             \& Centre for Operational Research, Management Sciences and Information Systems (CORMSIS)\\
             University of Southampton\\
             Building 54 Mathematical Sciences, SO17 1BJ Highfield Campus\\
              Tel.: +44(0)23 8059 3863\\
             % Fax: +123-45-678910\\
              \email{a.b.zemkoho@soton.ac.uk}           %  \\
%%             \emph{Present address:} of F. Author  %  if needed
}

\date{Received: date / Accepted: date}

\maketitle

\begin{abstract}
The aim of this paper is to present a set of Python-based tools to develop forecasts using time series data sets. The material is based on a four week course that the author has taught for seven years to students on operations research, management science, analytics, and statistics one-year MSc programmes. However, it can easily be adapted to various other audiences, including executive management or some undergraduate programmes. No particular knowledge of Python is required to use this material. Nevertheless, we assume a good level of familiarity with standard statistical forecasting methods such as exponential smoothing, AutoRegressive Integrated Moving Average (ARIMA), and regression-based techniques, which is required to deliver such a course. Access to relevant data, codes, and lecture notes, which serve as based for this material are made available (see \href{https://github.com/abzemkoho/forecasting}{github.com/abzemkoho/forecasting}) for anyone interested in teaching such a course or developing some familiarity with the mathematical background of relevant methods and tools.
\keywords{Python \and forecasting \and exponential smoothing \and ARIMA \and regression}
\subclass{90-04 \and  97U50 \and  97U70}
\end{abstract}

%\listoftables
%\vspace{5pt}

\section{Introduction}\label{intro}
According to Makridakis et al. \cite{MakridakisEtAl1998}, forecasting is the \emph{process of making predictions of the future based on past and present data and most commonly by analysis of trends, and usually needed to determine when an event will occur or a need arise, so that
appropriate actions can be taken}. Broadly speaking, forecasting approaches can be split into two categories: \emph{qualitative} and \emph{quantitative} forecasting methods, and we can even add a third one, that we label as \emph{semi-qualitative}, where a combination of both qualitative and quantitative methods can be employed to generate forecasts. Qualitative forecasting methods are often used in situations where historical data is not available. For more details on these concepts, interested readers are referred to the books \cite{HyndmanEtAl2018,MakridakisEtAl1998} and references therein. % Further details on qualitative forecasting methods can be found in Chapter
%3 of the book by Hyndman and Athanasopoulos (2014).

 Our focus in this paper is on quantitative methods, as we assume that historical \emph{times series data} (i.e., data from a unit (or a group of units) observed in several successive periods) is available for the variables of interest. Within quantitative methods, we also have a number of subcategories that can be broadly labelled as \emph{statistical methods}, which are at the foundation of the subject, and \emph{machine learning} ones, which have been developing rapidly in recent years; see, e.g., \cite{Hamzac,Deng,Robinson,Salaken,Voyant,Zhang,Zhang2} for a sample of applications and surveys on the subject. %These latter categories certainly overlap in some way, with techniques based regression analysis, for example [...].

The material to be presented in this paper is based on statistical forecasting methods; see, e.g., \cite{Adya,Chatfield,HyndmanEtAl2018,MakridakisEtAl1998,Sharda} for related details. Despite the fast development of machine learning techniques, they have been consistently shown through the last two \emph{M} competitions \cite{MakridakisEtAl2018,MakridakisEtAl2020} to generally be outperformed by statistical methods in terms of accuracy and computational requirements; these comparisons (see relevant details in the papers  \cite{MakridakisEtAl2018,MakridakisEtAl2020}) are done on more than 100 thousand practical data sets, related to a wide range of industries, based on the ForeDeCk database (\href{http://fsudataset.com/}{fsudataset.com}). Note that the M competition series (with M referring to Spyros Makridakis, one of the world leaders in the field) is a famous open competition, which can also be seen as a benchmarking exercise, where competitors evaluate and compare the performance of a wide range of forecasting methods on thousands of practical data sets.

The aim of this paper is to introduce the reader to existing Python tools that can be used to deliver a practical course on basic statistical forecasting methods; namely, we will focus on the exponential smoothing, AutoRegressive Integrated Moving Average (ARIMA), and regression-based methods, which are (or a combination of them) part of core techniques shown to have the best performance in the M competitions mentioned above.

\subsection{Background} The material presented in this paper is based on a course named \emph{Forecasting}, that the author has taught for the past seven years within the School of Mathematical Sciences at the University of Southampton, based in the United Kingdom. This  is an optional course, but which is very popular,  and is taken by students from the eight MSc programmes listed in Table \ref{table:1}, spanning both the School of Mathematical Sciences and the Southampton Business School.
%\begin{center}
%\begin{tabular}{| l| l| }
% \hline
%\textbf{School of Mathematical Sciences} & \textbf{Southampton Business School} \\
%  \hline
%  Data and Decision Analytics & Business Analytics and Management Science \\
% Operational Research & Business Analytics and Finance \\
%  Operational Research and Finance & Business Analytics and Finance \\
%   Statistics & Marketing Analytics\\
% \hline
%\end{tabular}
%\end{center}

\begin{table}[h!]
\centering
\begin{tabular}{||l l||}
 \hline
 \emph{School of Mathematical Sciences} & \emph{Southampton Business School}  \\ [0.5ex]
 \hline\hline
Operational Research and Finance & Business Analytics and Management Science \\
Data and Decision Analytics   & Logistics and Supply Chain Analytics \\
 Operational Research  & Business Analytics and Finance \\
   Statistics & Marketing Analytics\\[1ex]
   \hline
\end{tabular}
\caption{List of MSc programmes of origin of the students that usually take the forecasting course, which is the source of the material presented in this paper.}
\label{table:1}
\end{table}
${}$\\[-6ex]
%\begin{itemize}
%  \item Data and Decision Analytics;
%  \item Operational Research;
%  \item Operational Research and Finance;
%  \item Business Analytics and Management Science;
%  \item Business Analytics and Finance;
%  \item Logistics and Supply Chain Analytics;
%  \item Marketing Analytics,
%  \item Statistics.
%\end{itemize}
The course is very practical and hands on, designed to run for 16h across 4 weeks, with 2h of weekly lecture  and the remaining 2h dedicated to a  workshop/tutorial/computer lab, where the students are supported to go through the Python material to test and apply the methods on some practical data sets. The lectures focus on taking the students through the mathematical background of the methods that will be covered here \cite{Zemkoho2021}. During the computer labs, students are taken through the Python codes covered in this paper, which implement the methods that form the content of the lectures, and support them in using these methods to develop forecasts on practical data sets. Note that this course can easily be expanded to cover a few more weeks, as necessary, and the material can also be adapted to an undergraduate level for programmes around operations research, statistics, business analytics, and management science.

{ {It is important to mention that before the start of the course, a brief material with a basic introduction to Python is made available to the students, in order to bring them up to speed with some basic elements of Python, in case they have had no prior exposure to the language. This brief material essentially covers the relevant Python ecosystem discussed in Section \ref{Preliminaries on Python and data analysis} and an overview of the basic steps needed to get Python up and running on their personal computers or the university machines. Additionally, note that each of the weekly computer labs, which take place during the course, is an opportunity for the instructors to guide the students on how to use the different libraries needed to implement the mathematical concepts covered in the lecture of that week.}}

The author has taught the course over the last seven years,  first using Excel and relevant Visual Basic for Applications (VBA) codes to enhance some of the techniques. The transition to Python was done much recently, considering the demand both from industry and students, and also to keep up with the pace of  developments in data science more broadly. The motivation to prepare this paper came as a result of the transition from Excel to Python, as the author was unable to find a single book or resource relevant to prepare for a complete delivery of this course using Python.% where Python is used for all the methods presented here.

\subsection{Contribution of the paper and relevant literature}
The paper will mostly focus on the use of existing Python tools to generate forecasts, although a bit of the background on the mathematical concepts will be provided as necessary. Also, although prior knowledge of Python is not necessary, it will be assumed that the reader has some level of familiarity with methods involved in the corresponding mathematical material, as it would be required for anyone teaching such a course. The lecture notes \cite{Zemkoho2021} that form the material of the course discussed here are based on the books \cite{HyndmanEtAl2018,MakridakisEtAl1998}.

As for the Python material, we only found the book \cite{Brownlee2021} during the preparation of the first draft of the computing material, to be presented here, in 2019. While preparing this paper, we came across the two new books \cite{Korstanje2021,Lazzeri2021} on the use of Python to generate forecasts on time series data. There are two common denominators to  these three books; the first is that they are mostly geared towards machine learning-based techniques for time series forecasting, with the exception of ARIMA models, which are covered in detail. Secondly, they essentially focus on the use of Python tools to generate forecasts, and hence do not specifically pay attention to the mathematical background of the methods, which are the based on the corresponding Python forecasting tools.

Clearly, there are two differences between the content of this paper and what is covered in the books \cite{Brownlee2021,Korstanje2021,Lazzeri2021}. At first, considering the page limitation of an article such as this one, we also mostly only focus on the coding side of the methods; however, our presentation is essentially organized along the lines of the corresponding lecture notes \cite{Zemkoho2021}, which provide the necessary mathematical background to develop a deep understanding of all the methods covered in this paper. Secondly, unlike in these books, we focus our attention on statistical methods, which form the basis of most of the methods which are at the heart of the successful practical implementations in the context of the M competition series, as discussed at the beginning of this introduction.

%the book is largy focused on machine learning methods, and for the statistical techniques covered in this paper, only ARIMA approaches are covered but with codes based on [ML]? The drawback of these two books is that they mostly focus only on the use of existing Python tools to generate forecasts and the corresponding mathematical component, base of these methods, is generally very week.
%
%The book by Jason Brownlee [...] can be important for some aspects of the material in this paper, in particular, for Sections ... Also, while preparing this paper, the author came accross the book [..] just published last month (July 2021). However, more broadly, the author is not aware of any resource that cover the topics to be covered here, in such a concise way.

It is also important to mention that our philosophy in the preparation and delivery of the course discussed in this paper is inspired in part by the book \cite{HyndmanEtAl2018}; that is, giving the reader a balanced mathematical background of the forecasting methods, while accompanying them with relevant practical software tool to use these methods on practical data sets. However, the fundamental difference is that \cite{HyndmanEtAl2018} uses R while we use Python.

The lecture notes on which this course is based (i.e., \cite{Zemkoho2021}), as well as all the corresponding codes presented here, can be accessed online via the following link: \href{https://github.com/abzemkoho/forecasting}{github.com/abzemkoho/forecasting}.

%\subsection{Relevant literature}
%----------------------

%The forecasting methods that form the base for the material presented in this paper are based on the book [...] The module uses Makridakis, S., Wheelwright, S.C. and Hyndman, R.J. 1998, Forecasting:
%Methods and Applications 3rd Ed., New York: Wiley as text book. Most of the material for the
%lecture notes is extracted from there. Also, most of the data sets used in the demonstrations are
%drawn from this book. (The full set from the book can also be downloaded under Course Content,
%if desired.) Hyndman, R.J. and Athanasopoulos, G. 2014, Forecasting: principles and practice,
%Ortexts.com, although slightly lighter, has recently superseded the book by Makridakis et al. Hence,
%some of the material of these notes has also been drawn from there. An additional advantage of
%the book by Hyndman and Athanasopoulos (2014) is that it is freely accessible to read online (at
%{https://www.otexts.org/fpp/}). This book also has a hands on approach relevant to this paper. However, the codes there are based on the R package.

\subsection{Outline of the paper}
We start the next section with an overview of the main Python packages needed to work with the tools that we will go through in this paper. Subsequently, we present tools that can be used for a basic data analysis (i.e., time, seasonal, and scatter plots, as well as correlation analysis, just to mention a few) before the start of any forecasting task based on the methods  covered in this paper. Section \ref{Exponential smoothing methods0} is devoted to exponential smoothing methods, which are very efficient on time series that involve trends and/or seasonality. Section \ref{ARIMA methods} covers autoregressive integrated moving average (ARIMA) methods; and finally, section \ref{Regression-based forecasting method} presents tools for regression analysis and how they can be used for forecasting. Note that the exponential smoothing and ARIMA methods are blackbox techniques, as are they are built under the assumption that historical patterns in the time series will keep repeating themselves in the future. However, regression-based approaches assume that the behaviour of the time series of interest (dependent variable) is influenced by other variables (independent variables), and this is explored through linear regression to possibly build more accurate forecasts. % where it assumed that .... that ..... However, the methods ....

%\begin{itemize}
%  \item Say a bit more about each method
%\end{itemize}

\section{Preliminaries on Python and data analysis}\label{Preliminaries on Python and data analysis}
%\begin{itemize}
%  \item It seems like the books heavily rely on the scikit-learn library
%  \item we mostly use Pandas and Statsmodels?
%  \item Two { {SciPy}} libraries provide a foundation for most others: they are \textbf{NumPy} (for providing efficient array operations) and \textbf{Matplotlib} (for plotting data)
%  \item There are three higher-level SciPy libraries that provide the key
%features for time series forecasting in Python: they are \textbf{Pandas, Statsmodels, and scikit-learn}
%for data handling, time series modeling, and machine learning respectively. Let's take a closer
%look at each in turn.
%\end{itemize}
\subsection{The necessary Python ecosystem}
No prior knowledge of Python is required to use the material in this paper. However, we assume that the reader/instructor who wants to use the tools presented here has Python up and running on their device (desktop, laptop, etc.) The codes and corresponding results  are based the use of Python under \texttt{Anaconda 3} with \texttt{Spyder 3.6} as editor, all running on Windows 10 Enterprise (processor: Intel(R) Core(TM) i5-6300U CPU \@ 2.40 GHz). The advantage of using Anaconda is that it installs Python with many important packages that are useful for time series analysis of the type covered in this paper. This therefore helps in part to reduce dependency issues between various packages used, and hence ensure that key packages are set to work nicely together. Nevertheless, all the codes presented here should be able to work smoothly on most platforms running a version 3 of Python (see \href{https://www.python.org/}{python.org}). The main packages needed are as follows:
\begin{itemize}
  \item SciPy;
  \item NumPy;
  \item Matplotlib;
  \item Pandas;
  \item Statsmodels.
\end{itemize}
As an \emph{ecosystem of open-source software for mathematics, science, and engineering}, SciPy (\href{https://scipy.org/}{scipy.org}) includes the NumPy library (\href{https://numpy.org/}{numpy.org}) for  multi-dimensional array operations and  the Matplotlib library (\href{https://matplotlib.org/}{matplotlib.org}) specifically designed for plotting.

Pandas (\href{https://pandas.pydata.org/}{pandas.pydata.org}) is an open source data analysis and manipulation tool.  It is important to mention here that all the data sets used for our illustrations are stored in Excel spreadsheets. Hence, we use pandas function named \texttt{read\_excel} in almost all our codes, to read data from Excel worksheets.
Finally, the statsmodels library (\href{https://www.statsmodels.org/stable/index.html#}{statsmodels.org}) is at the heart of most of the data analysis and forecasting methods presented in this paper, as it includes packages to generate forecasts using exponential smoothing, ARIMA, and regression-based methods. Occasionally, we could also use \texttt{scikit-learn} for a few tasks, such as generating some error measures.

It is also important to mention that the above mentioned libraries all work together, with NumPy and Matplotlib building on SciPy, while Statsmodels is built on top of NumPy and SciPy, integrating with Pandas for data handling, as already mentioned.

%We should also make sure that pandas, numpy, etc. are installed ... It could be that pandas, numpy, etc, are part of most standard installations of a Python 3x. As for the stats... package, it might be necessary to install separately.
%
%We will also assume that all the data sets to be used in our analysis are stored in Excel spreadsheets; hence, the reason why most of our codes will involve the following command lines:\\[-1ex]
%
%\begin{lstlisting}[language=Python]
%from pandas import read_excel
%series = read_excel('Electricity.xls', sheetname='SeasData', header=0, index_col=0, parse_dates=True, squeeze=True)
%\end{lstlisting}
%(say a few words here about the methods)
%\\[2ex]
%
%\noindent \textbf{Library: Statsmodels} -- used for the following stuff
%\begin{enumerate}
%  \item Section 2 (preliminary data analysis)
%  \begin{itemize}
%    \item to assess seasonality (see section 2)
%    \item to decompose a time series (see section 2)
%  \end{itemize}
%  \item Section 3 (exponential smoothing methods)
%  \begin{itemize}
%    \item Base for all the exponential smoothing methods (see section 3)
%    \item (SimpleExpSmoothing, Holt, ExponentialSmoothing)
%  \end{itemize}
%  \item Section 4 (regression analysis)
%  \begin{itemize}
%    \item ols function needed everywhere
%  \end{itemize}
%  \item Section 5 (ARIMA)
%  \begin{itemize}
%    \item ACF and PACF plots
%  \end{itemize}
%\end{enumerate}

\subsection{Basic data analysis tools}\label{Basic data analysis tools}
\begin{figure}[htp]
\hfill
\subfigure[Clay bricks sale]{\includegraphics[width=4.7cm]{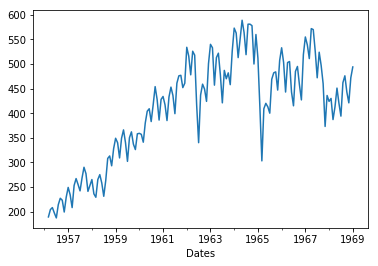}}
\hfill
\subfigure[Electricity demand]{\includegraphics[width=5cm]{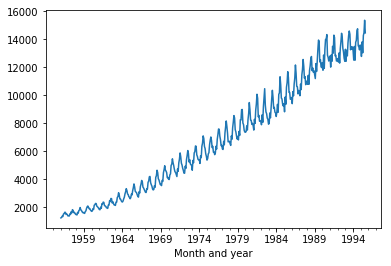}}
\hfill
\subfigure[USA treasury bills]{\includegraphics[width=4.7cm]{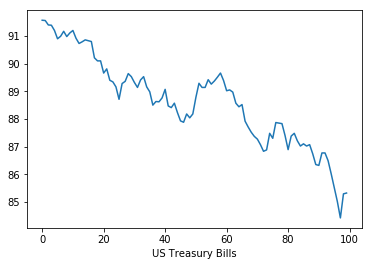}}
\hfill
\caption{Time plots for three data sets demonstrating different patterns; as for all the data sets used in this paper, they are drawn from material related to the book \cite{MakridakisEtAl1998} and will be made available in the supplementary files.}\label{graph0}
\end{figure}

In this subsection, we discuss the following five key topics, which are crucial in the preliminary analysis of time series data sets:
\begin{itemize}
  \item time plots;
  \item adjustments;
  \item decompositions;
  \item correlation analysis;
  \item autocorrelation function.
\end{itemize}
A time plot is simply a two-dimensional graphical representation of a time series. A time plot is typically the starting point of any forecasting task, as it enables a first dive into the data set, to assess, for example, whether any errors or particular patterns are present in it. To build a time plot, \texttt{matplotlib} can be used after the data has been loaded using the \texttt{read\_excel} function from pandas, as described above. Listing \ref{TimePlot} provides the code for the time plot in Figure \ref{graph0} (a); the remaining ones are obtained just by replacing the clay bricks data set by the corresponding ones. Note that \texttt{read\_excel} includes a number of options to specify the sheet and other information to be returned for use by the \texttt{series.plot()} command, which generates the plot.

%\begin{figure}
%  \centering
%\includegraphics[width=\textwidth]{Demo11Fig1}
%\caption{Overall, there is a global trend with significant trend up to 1965, where we start having some occasional large fluctuations which
%are difficult to explain, and hence predict, without knowing the underlying causes, but could be a reflection of a cyclical behaviour.}\label{graph0}
%\end{figure}

Obviously, Figures \ref{graph0} (a) and (b) show a globally increasing trend, while the trend in (c) is generally decreasing. However, in the clay bricks sale, there is overall a global trend with steady increase over time up to 1965, where we start having some occasional large fluctuations which
are difficult to explain, and hence, hard to predict without knowing the underlying causes. But this  could be a reflection of a cyclical trend behaviour in the  clay bricks sale time series.

Clearly, trend and cyclical behaviour can be easily identified from a graph. However, unlike trends and cyclical patterns, seasonality can be trickier to observe. Considering the important role that the identification of seasonality plays in developing/selecting some forecasting methods (e.g., exponential smoothing and ARIMA methods, as it will be clear in the following sections), we need to pay some attention on how to assess it. There are various ways to assess whether a time series has seasonality, including zooming out specific chunks of the corresponding time plots. Also, a time plot can sometimes already give an initial indication on the presence of seasonality in a time series; for example, intuitively, Figure \ref{graph0} (b) already suggest that we might be having peaks and troughs occurring at regular intervals. But some further steps need to be taken to check this.

%(see, e.g., \href{https://www.dataquest.io/blog/tutorial-time-series-analysis-with-pandas}{https://www.dataquest.io/blog/tutorial-time-series-analysis-with-pandas} ).
In this paper, we are going to mainly use the \emph{seasonal plots} and the concept of \emph{autocorrelation function} (ACF) to decide whether a time series is seasonal or not. The ACF will be defined at the end of this section. Before that, we start with the seasonal plots, which correspond to a superposition of time plots over a succession of limited time periods (e.g., 12 months in the context of monthly observations, which is what we have for most of the data sets used in our illustrations). Listing \ref{SeasonalPlot} provides a code that can be used to build seasonal plots after having organized our data in months for over a few years. %, as all our data sets are made of monthly observations. Hence, we do the superpositions over regular intervals of 12 months.
%We can clearly see that the time series is cut into regular periods and the time plots of each period are overlaid on top of
%one another. It is an effective means of demonstrating seasonality, for example to a client who is
%not overly technically minded.
\begin{figure}[htp]
\hfill
\subfigure[Clay bricks sale]{\includegraphics[width=4.7cm]{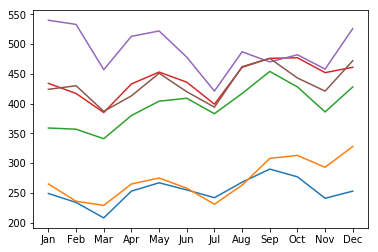}}
\hfill
\subfigure[Electricity demand]{\includegraphics[width=5cm]{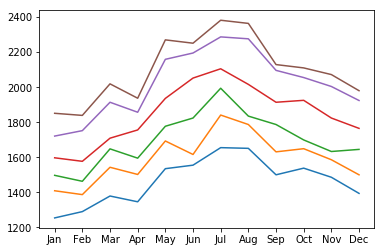}}
\hfill
\subfigure[USA treasury bills]{\includegraphics[width=4.7cm]{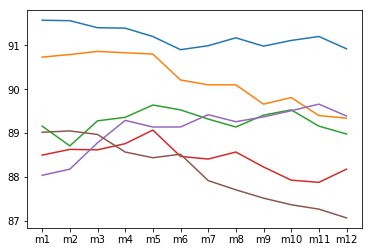}}
\hfill
\caption{Seasonal plots for the time series plotted in Figure \ref{graph0}; the level variation in each of the graphs is further evidence of the increasing or decreasing trends that clearly appear in the time plots.}\label{graph1}
\end{figure}
%\begin{figure}
%  \centering
%\includegraphics[width=\textwidth]{Demo12Fig1}
% \caption{Overall, there is a global trend with significant trend up to 1965, where we start having some occasional large fluctuations which
%are difficult to explain, and hence predict, without knowing the underlying causes, but could be a reflection of a cyclical behaviour.}\label{graph1}
%\end{figure}

Clearly, there is an indication from Figure \ref{graph1} that the clay bricks and electricity data that may have seasonality, while it is unlikely to be the case for the treasury bills data. From the time plots in Figure \ref{graph0}, an initial guess could have already been made about the electricity data, but maybe not necessarily for the clay bricks data. At the end of this section, we will see how the ACF plots can help to further confirm seasonality identified here.

%If seasonal plots lead you to think that seasonality may exist, then further technical evidence can be
%provided by autocovariance and autocorrelations. These concepts are introduced in the next section.
%Before going to that, we provide another class of plots that is useful in analyzing bivariate-type
%data sets (i.e., data sets based on two variables).
%
%For time and seasonal plots, we essentially need the \texttt{plot} function from \texttt{matplotlib}, which can be applied as
%\begin{lstlisting}[language=Python]
%series.plot()
%\end{lstlisting}
%with series representing the time series data set to be represented. The difference in the way both types of plots are used is that the time plot represents the data set completely, while in the case of seasonality, we just have a superposition of time plots over a succession of limited time periods (e.g., 12 months in the context of monthly observations). Appendices ... and ... provide complete examples illustrating how time and seasonal plots can be built.
%
%--------------------------------

Besides the different patterns that can be assessed using time plots, they can also enable an assessment of the need for adjustments (e.g., mathematical transformations or calendar adjustments). Ideally, the role of a mathematical transformation is to attempt to stabilize variance in a time series, where rapid changes in some parts of a time plot can affect the ability of a forecasting method to generate accurate results. For instance, the power (including the square root, as
%\begin{figure}[htp]
%\hfill
%\subfigure[Clay bricks sale]{\includegraphics[width=3.1cm]{Demo12Fig11}}
%\hfill
%\subfigure[Electricity demand]{\includegraphics[width=3.1cm]{Demo12Fig12}}
%\hfill
%\subfigure[USA treasury]{\includegraphics[width=3.1cm]{Demo12Fig13}}
%\hfill
%\subfigure[USA treasury]{\includegraphics[width=3.1cm]{Demo12Fig13}}
%\hfill
%\caption{Overall, there is a global trend with significant trend up to 1965, where we start having some occasional large fluctuations which
%are difficult to explain, and hence predict, without knowing the underlying causes, but could be a reflection of a cyclical behaviour.}\label{graph2}
%\end{figure}
 a special case) and log transformations are the most commonly used transformations in the literature; the square root can help, in the case where the time series has the shape of a second order quadratic function, to promote a ``linear'' shape, which can improve the predictability capacity of some forecasting methods. On the other hand, the log (of course, applicable only for positive time series) has an additional advantage, in terms of its interpretability. For more details on these transformations and many other adjustments, which can positively impact the forecasting ability of some methods, see \cite[Chapter 3]{HyndmanEtAl2018}. Listings \ref{LogTransform}, \ref{SqrtTransform}, and \ref{CalendarAdjustment} provide appropriate codes to generate a log, square root, and calendar adjustments, respectively. %, for a milk data sets.
The code in Listing \ref{CalendarAdjustment} runs on a special data set, where a calendar adjustment can be useful, as in the milk production of a cow, the difference in the observations from one month to the other can essentially be due to the number of days in months. Hence, the calendar adjustment can help to remove such a calendar effect before any further analysis of this time series.

\begin{figure}[htp]
  \begin{center}
    \includegraphics[width=0.50\textwidth]{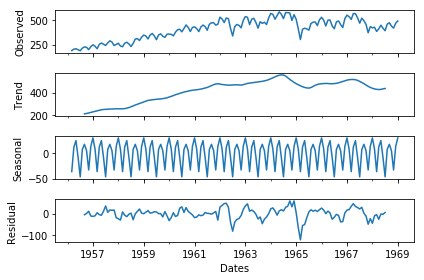}
  \end{center}
  \caption{Additive decomposition graphs for the clay bricks sale time series}\label{Decompose11}
\end{figure}

For a given time series $\{Y_t\}_t$, it is sometimes important to look for ways to split it by means of a decomposition function $f$ in such way that
\begin{equation}\label{Decompo}
  Y_t = f(T_t, S_t, E_t),
\end{equation}
where for a given $t$, $T_t$ and $S_t$ denote the trend-cycle and seasonal components, respectively, and $E_t$ corresponds to the error that results from such a decomposition. Decompositions are useful in developing a better understanding of the constituting patterns in a time series, but not necessarily for generating forecasts. Standard selections for a decomposition function are $ f(T_t, S_t, E_t):= T_t + S_t + E_t$ (additive decomposition) and  $ f(T_t, S_t, E_t):= T_t \times S_t \times E_t$ (multiplicative decomposition). The \texttt{statsmodels} function \texttt{seasonal\_decompose} can be used to generate these decompositions, with the option ``model'' suitable for indicating the nature of the decomposition (i.e., additive or multiplicative); see Listing \ref{DecompositionAdditive} for an additive decomposition code
 (used to generate Figure \ref{Decompose11}, for illustrative purpose) and Listing \ref{DecompositionMultiplicative} for a multiplicative one.

%
%\begin{figure}[htp]
%\hfill
%\subfigure[Clay bricks sale]{\includegraphics[width=4.7cm]{Decompose11}}
%\hfill
%\subfigure[Electricity demand]{\includegraphics[width=5cm]{Demo12Fig12}}
%\hfill
%\subfigure[USA treasury bills]{\includegraphics[width=4.7cm]{Demo12Fig13}}
%\hfill
%\caption{Seasonal plots for the time series plotted in Figure \ref{graph0}; the level variation in each of the graphs is further evidence of the increasing or decreasing trends that clearly appear in the time plots.}\label{graph1}
%\end{figure}
It is important to note that in terms of the background algorithm on how a decomposition is computed, one usually starts with the trend estimation, and then, depending on the nature of $f$ (\ref{Decompo}), the seasonal component is estimated; interested readers are refereed to the lecture notes associated to this material \cite[Section 2]{Zemkoho2021} and references therein.

%%--------------------------------
%%
%%The calendar adjustment plays a similar role as the transformations, with the difference that ....
%%
%%To plot transformations and calendar adjustments, the corresponding transformations first need to be applied to the original time series before the plots are constructed; see Listings ... and ... for the examples of illustrating codes.
%
%
%----------

Correlation analysis comes into play when we want to explore relationships between variables in cross-sectional data. There are at least two possible tools to assess correlation between variables. Namely, scatter plots and correlation values; both concepts are strongly related in the sense that the scatter plot provides a graphical representation that can demonstrate how strong the relationship between two variables is, while the correlation is a numerical value materializing the strength level of such a relationship. As, example to illustrate these two concepts, consider a data set made of a variety of used cars and their price (based on their mileage). For instance, we might want to forecast (price) against one possible explanatory variables (mileage, here). Running the code in Listing \ref{correlation} clearly shows that the price of a car decreases as the mileage increases.  Each point on the graph represents one specific vehicle.

\begin{figure}[htp]
\hfill
\subfigure{\includegraphics[width=6.5cm]{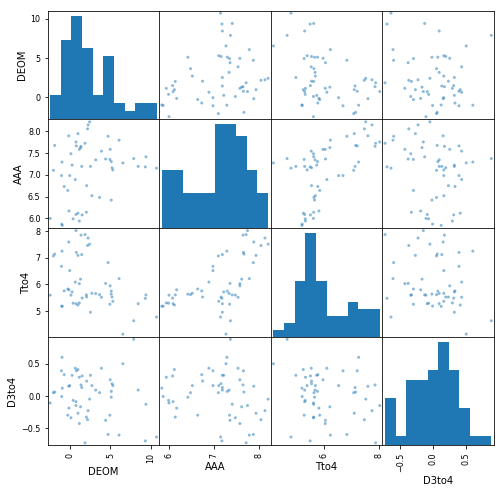}}
\hfill
\subfigure{\includegraphics[width=8cm]{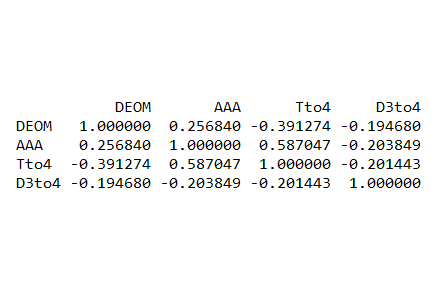}}
\hfill
\caption{Left, we have the matrix of scatter plots for four times series labelled as DEOM, AAA, Tto4, and D3to4. On the right, we have the correlation matrix, which gives the correlation value that reflects the relationship in each pair in these four data sets. As it can be seen in the scatter plots, the strongest correlation is between AAA and Dto4, as confirmed by the correlation value, which is strictly larger than $0.50$.}\label{graphCorrelation}
\end{figure}

A scatter plot helps us to visualize the relationship and suggests
that if one wants to forecast the price of used car, a suitable model should include mileage as an explanatory variable. In Listing \ref{correlation}, the scatter plot function \texttt{scatter} function from \texttt{matplotlib} is applied with arguments being the mileage and price as separate entries. Note that \texttt{pandas} also has the function \texttt{scatter\_matrix}, which can generate scatter plots for many variables in one go; this could be particularly important in Section \ref{Regression-based forecasting method} when studying the regression approach to forecasting. Figure \ref{graphCorrelation}, for example, generated by the code Listing \ref{CorrelationMatrix}, shows scatter plots in a matrix form for four time series.

The \emph{correlation} is a statistic corresponding to a number between -1 and 1 to measure the level of the linear relationship for bivariate data (i.e.,  when there are two variables). The \texttt{corrcoef} function from \texttt{numpy}, see Listing \ref{correlation}, calculates the correlation between the mileage and prices of the cars, as discussed above. Note that in principle, \texttt{corrcoef} is generated as a symmetric matrix; hence the use of \texttt{correlval}[1,0] to extract the necessary value. In a situation where one is interested in evaluating the relationships between various pairs of variables, the correlation matrix enables the calculation of these values in one go, as discussed above in the context of scatter plots, as illustrated in the left-hand-side of Figure \ref{graphCorrelation}; the corresponding correlation values are generated with the function \texttt{corr} from \texttt{pandas}; see the table in the right-hand-side of Figure \ref{graphCorrelation} for an illustration with four time series.

%Hence, \texttt{corrcoef} can be applied to the case where one has multiple variables (see regression chapter). There are also other functions that can be used? See regression and other places where correlation are calculated.
%
%I think it would be better to have a separate small section on correlation, where correlation between 2 or more variables will be analyzed from the perspective of (1) scatter plots and (2) correlation matrix; see notes and Demo 3.1-2 (Correlation matrix), corr (), directly from pandas is used. It seems corroef is specifically from numpy.
%
%Could we also do same for seasonality? i.e., present it with both seasonal plots and autocorrelation function?
%
%For a single time series, the concept of correlation can also be used to look for specific patterns. As example, in the example in Listing \ref{autocorrelation1} demonstrates that there is seasonality in the beer sales in Australia.
%
%
%----
%(Calculating/plotting/representing the autocorrelation)
%----
%

\begin{figure}[htp]
\hfill
\subfigure{\includegraphics[width=6.5cm]{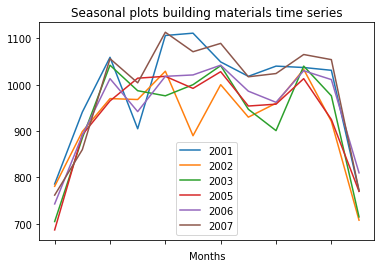}}
\hfill
\subfigure{\includegraphics[width=6.5cm]{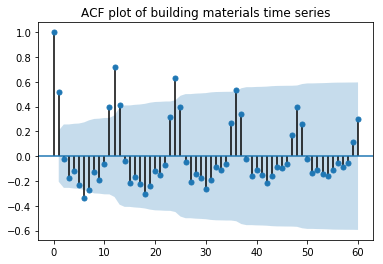}}
\hfill
\caption{Left, we have the seasonal plots for most of the years involved in the times series. On the right-hand-side, we have the ACF plot over 60 time lags.}\label{graphAutoCorrelation}
\end{figure}
For a given time series $Y_t$, the concept of correlation can be extended to the time lags $Y_t$ and $Y_{t-k}$ of this same series. Hence, such a correlation is called \emph{autocorrelation}. The autocorrelation is used to measure the degree of correlation between different time lags in a time series. The autocorrelation function (ACF) is crucial in assessing many properties in statistics, including seasonality, white noise, and stationarity. In this section, we limit ourselves to the use of the ACF in assessing seasonality. For its use in assessing white noise and stationarity, see Sections \ref{Exponential smoothing methods0} and \ref{ARIMA methods}, respectively.

The function \texttt{plot\_acf} from \texttt{statsmodels} can be used to generate ACF plots. As illustration, we use the code in Listing \ref{Autocorrelation} to generate the ACF plot for some building material data from Australia in Figure \ref{graphAutoCorrelation}. It can be seen that the seasonal plots (left-hand-side) might be slightly unclear. But the ACF plot (right-hand-side) strongly confirms the presence of seasonality, as peaks and troughs are approximately occurring at regular interval; namely, at every 12 time lags, as the data is made of monthly observations. Note that \texttt{autocorrelation\_plot} from \texttt{pandas} can also be used to produce ACF plots, but instead with a curve shape; this can be seen by running the code in Listing \ref{Autocorrelation}, where the relevant command is also included. %. The function \texttt{plot\_acf} could also be used to generate the ACF graph in a histogram format. It is clear that both graphs have the same shapes.

%%%%%%%%%%%%%%%%%%%%%%%%%%%%%%%%%%%%%%%%%%%%%%%%%%%%%%%%%%%%%%%%%%%%%%%%%%%%%%%%%%
%%%%%%%%%%%%%%%%%%%%%%%%%%%%%%%%%%%%%%%%%%%%%%%%%%%%%%%%%%%%%%%%%%%%%%%%%%%%%%%%%%
%%%%%%%%%%%%%%%%%%%%%%%%%%%%%%%%%%%%%%%%%%%%%%%%%%%%%%%%%%%%%%%%%%%%%%%%%%%%%%%%%%
%%%%%%%%%%%%%%%%%%%%%%%%%%%%%%%%%%%%%%%%%%%%%%%%%%%%%%%%%%%%%%%%%%%%%%%%%%%%%%%%%%
%%%%%%%%%%%%%%%%%%%%%%%%%%%%%%%%%%%%%%%%%%%%%%%%%%%%%%%%%%%%%%%%%%%%%%%%%%%%%%%%%%
%%%%%%%%%%%%%%%%%%%%%%%%%%%%%%%%%%%%%%%%%%%%%%%%%%%%%%%%%%%%%%%%%%%%%%%%%%%%%%%%%%
%%%%%%%%%%%%%%%%%%%%%%%%%%%%%%%%%%%%%%%%%%%%%%%%%%%%%%%%%%%%%%%%%%%%%%%%%%%%%%%%%%
%%%%%%%%%%%%%%%%%%%%%%%%%%%%%%%%%%%%%%%%%%%%%%%%%%%%%%%%%%%%%%%%%%%%%%%%%%%%%%%%%%
%%%%%%%%%%%%%%%%%%%%%%%%%%%%%%%%%%%%%%%%%%%%%%%%%%%%%%%%%%%%%%%%%%%%%%%%%%%%%%%%%%
%%%%%%%%%%%%%%%%%%%%%%%%%%%%%%%%%%%%%%%%%%%%%%%%%%%%%%%%%%%%%%%%%%%%%%%%%%%%%%%%%%
\section{Exponential smoothing methods}\label{Exponential smoothing methods0}
Considering the importance of accuracy in forecasting, we start this section by discussion a few tools that can be used to assess the accuracy of a method; namely, we discuss error measures, the concept of white noise via the ACF function, and confidence intervals. After this, in Subsection \ref{Exponential smoothing methods}, we introduce the exponential smoothing methods, which, together with the ARIMA methods represent the most widely used forecasting techniques in practice.
\subsection{Accuracy measures}
As accuracy is the first main concern when forecasting, we start here by discussing how some standard error measures, i.e., the  mean error (ME), mean absolute error (MAE), mean square error (MSE), percentage error (PE), mean percentage error (MPE), and the mean absolute percentage error (MAPE), can be computed using Python. To proceed, it is crucial to recall that an error measure on its own does not mean much, but rather, it can only make sense in a comparison setting of 2 or more methods. Hence, we introduce two na\"{\i}ve forecasting methods to illustrate how these error measures can be used in practice. We begin with a na\"{\i}ve forecasting, labelled as NF1, which assumes that for a times series $\{Y_t\}$, the forecast at time point $t+1$ is obtained as $F_{t+1} = Y_t$. Next, we consider a second na\"{\i}ve forecasting  method labelled as NF2:
\[
F_{t+1} = Y_t - S_t + S_{(t-12) + 1} \;\; \mbox{ with } \;\; S_t = \frac{1}{m+1}(mS_{t-12} + Y_t),
\]
where $S_t=Y_t$ for $t=1, \ldots, 12$ and with $m$ is the number of complete years of data available; for the initialization of the method, we set $F_{t+1}=Y_t$ for $t=1, \ldots, 12$.
\begin{table}[h!]
\centering
\begin{tabular}{l c r}
\begin{minipage}[b]{0.28\linewidth}
\includegraphics[width=\textwidth]{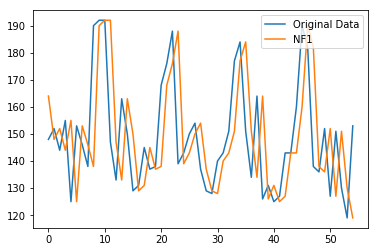}
%${}$\\[2ex]
\end{minipage} &
\begin{minipage}[b]{0.28\linewidth}
\includegraphics[width=\textwidth]{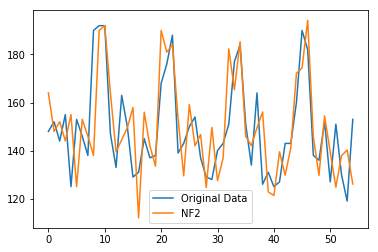}
%${}$\\[2ex]
\end{minipage}&
\begin{minipage}[b]{0.4\linewidth}
\includegraphics[width=\textwidth]{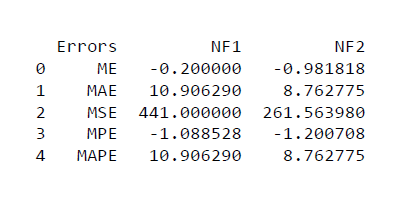}
%${}$\\[2ex]
\end{minipage}
\end{tabular}
\caption{The results from NF1 and NF2 can be seen in the first and second graphs, respectively. As for the corresponding error measures, see the table in right-hand-side.}
\label{table:3}
\end{table}

The code in Listing \ref{ErrorMeasures} generates the results in Table \ref{table:3}, which show both the NF1 and NF2 forecast plots, as well as the corresponding error measures stated above.  Note that the ME and MPE are not to be taken very seriously as their values essentially reflect the fact that positive and negative values just cancel each other throughout the range. Clearly, NF2 outperforms NF1 on almost all the measures, especially, on the positive ones (MAE, MSE, and MAPE), which are more meaningful. This is not surprising, considering the fact NF2 contains more structure capturing the nature of the data set much better than NF1, which is essentially a one step translation of the original data set. Similar comparisons can be done for any two or more forecasting methods.
\begin{figure}[htp]
\hfill
\subfigure[ACF plot of original data]{\includegraphics[width=4.7cm]{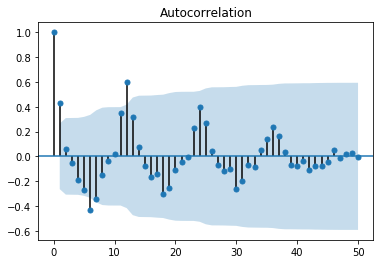}}
\hfill
\subfigure[ACF plot error from NF1]{\includegraphics[width=5cm]{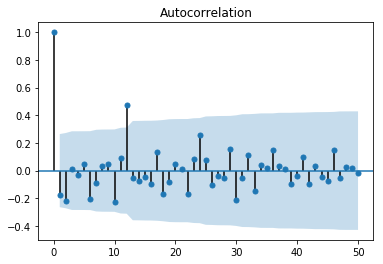}}
\hfill
\subfigure[ACF plot error from NF2]{\includegraphics[width=4.7cm]{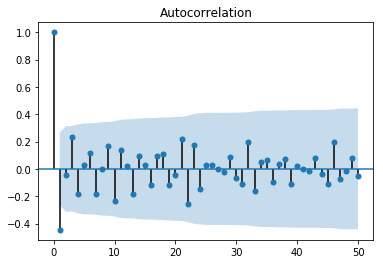}}
\hfill
\caption{The original time series is seasonal as it can be seen in (a) and this pattern is preserved in the errors resulting from NF1, with the spikes after each 12th time lag.}\label{graphACFerrors}
\end{figure}

Another tool to assess the accuracy of a forecast method is the ACF of the errors. Basically, the expectation is that if the results of a forecasting methods are reasonably accurate, the time plot of the errors, seen as a time series, should be purely random. Therefore, no patterns from the original data should be preserved in the errors/residuals. Using the corresponding code in Listing \ref{ACFErrors} on the data used for Figure \ref{table:3}, we get the graphs in Figure \ref{graphACFerrors}, which clearly show that the forecasts from NF1 preserve seasonality from the original time series, with the large spikes appearing after every 12th time lag. Such a pattern is not clearly obvious for NF2.

\begin{figure}[htp]
\hfill
\subfigure{\includegraphics[width=6.5cm]{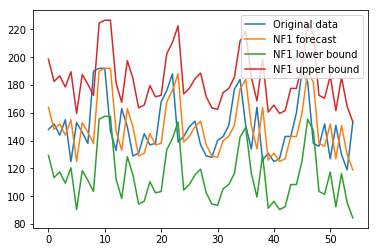}}
\hfill
\subfigure{\includegraphics[width=6.5cm]{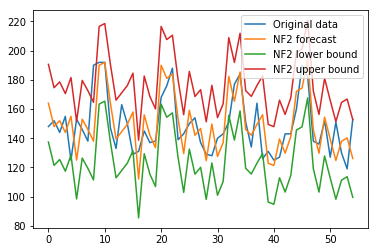}}
\hfill
\caption{The confidence intervals here are obtained with the formula $F_t \pm z\sqrt{MSE}$ with $z$ being the parameter ensuring that the $90\%$ chance that the forecasts would be between the lower and upper bounds provided.}\label{graphConfidenceInterval}
\end{figure}

Finally, providing the confidence interval for a forecast can help to a decision-makers in building their management perspectives.  Let $F_{t+1}$ be the forecast from a given method, then, the corresponding lower and upper bounds can be obtained as
\[
{LF}_{t+1}:= F_{t+1} - z\sqrt{\mbox{MSE}} \;\; \mbox{ and }\;\; {UF}_{t+1}:= F_{t+1} + z\sqrt{\mbox{MSE}},
\]
respectively, where MSE represents the mean square error over a suitable range of the data, while $z$ is a quantile of the normal distribution, which is a conventional number that determine the level of confidence of the corresponding interval. Standard values commonly used in practice for $z$ can be seen in Section 2 of \cite{Zemkoho2021}. Figure \ref{graphConfidenceInterval}, generated with the code in Listing \ref{ConfidenceInterval}, provides the confidence intervals for the data and corresponding NF1 and NF2-based results.

\subsection{Exponential smoothing methods}\label{Exponential smoothing methods}
There are four main types of exponential smoothing methods, which can be applied based on characteristics of our time series and sometimes also considering our intended purpose. Before diving into these methods, it is important to mention that all the related Python tools that we are going to describe here are from  the \texttt{statsmodels} library. The first and simplest such method is the so-called \emph{single exponential smoothing} (SES) method. The SES is usually applied only on time series that do not exhibit any specific pattern and can only produce one step ahead forecast.

To set the stage for the general process of all the forecasting methods that we are going to present in this paper, we are going to provide a brief overview of the mathematical background of the SES method. To proceed, let us assume that we are given a time series $Y_1$, \ldots, $Y_t$, where data is available from time point $T=1$ up to $T=t$. Then, the forecast for this time series at time point $T=t+1$ using the SES method can be calculated as
\begin{equation}\label{SESformula}
  F_{t+1} =(1-\alpha )^{t} F_{1} +\alpha \sum _{j=0}^{t-1}(1-\alpha )^{j} Y_{t-j},
\end{equation}
 where the parameter $\alpha\in [0, \; 1]$. There are various ways to initialize the method; one possibility is to select $F_1=Y_1$.
 The first key observation that can be made on the formula (\ref{SESformula}), and which justifies the name of this class of methods, is the fact that if one looks carefully at the factor $(1-\alpha)$, we will observe that it decays exponentially as the power $j$ increases. More interestingly, by the nature of the expression, this increase is associated with the decrease of the indices of $Y_j$. Hence, this means that the value of $F_{t+1}$ relies heavily on more recent values of the time series $Y_1$, \ldots, $Y_t$. This is one of the particular characteristics of any exponential smoothing methods.

\begin{figure}[htp]
\hfill
\subfigure{\includegraphics[width=5.5cm]{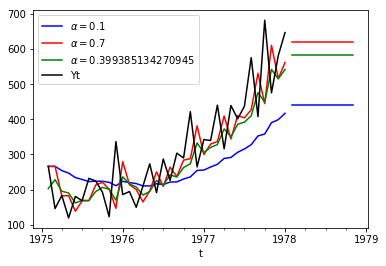}}
\hfill
\subfigure{\includegraphics[width=9.5cm]{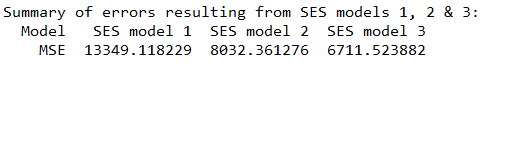}}
\hfill
\caption{On the left, we have the forecast plots for different values of the parameter $\alpha$, with the 3rd being the optimal one. The table on the right provides values of the MSE for each value of the parameter.}\label{graphSES}
\end{figure}

 Additionally, being able to optimally select the value of the parameter $\alpha$ is critical for the performance of the method. The strategy commonly used in this case is the least square optimization approach to select its best value. It corresponds to minimize the MSE
 \begin{equation}\label{SESoptimization}
 \min~\frac{1}{t}\sum^t_{j=1} e^2_j \;\; := \;\; \frac{1}{t}\sum^t_{j=1} \left(F_j - Y_j\right)^2 \;\;\mbox{ s.t. } \;\; \alpha\in [0, \; 1],
 \end{equation}
where $F_j$, $1=1, \ldots, t$ is obtained from (\ref{SESformula}).
The function from \texttt{statsmodels} to generate forecasts using the SES method is \texttt{SimpleExpSmoothing}. Applying it, as in Listing \ref{SimpleExpSmoothing}, to an example of data set generates the results in Figure \ref{graphSES}, where we can clearly see that the models based on manually selected values $\alpha=0.1$ and  $\alpha=0.7$ for models 1 and 2, respectively, are not as good as the 3rd model, which is obtained by solving the optimization problem in (\ref{SESoptimization}); clearly, the MSE values in the right-hand-side table in Figure \ref{graphSES} confirm that the parameter $\alpha$ obtained from the optimization problem (\ref{SESoptimization}) incurs the smallest error.

The core part of the code in Listing \ref{SimpleExpSmoothing} is \texttt{SimpleExpSmoothing} imported from the subpackage \texttt{statsmodels.tsa.api} of \texttt{statsmodels}. Everything else is essentially the selection of the parameter $\alpha$ for models 1 and 2 using the option \texttt{smoothing\_level}; obviously, in the case where $\alpha$ is optimally selected, the default selection of the \texttt{optimized} option is set at \texttt{True}. Also, recall that as the SES method can only produce a one-step forecast, the number 10 that appears in the command \texttt{fit2.forecast(10).rename(r'$\alpha=0.7$')}, in the 2nd model, for example, sets the number of times that the single value of $F_{t+1}$ is going to be repeated. This is essentially to clearly visualize the result in the graph; however, it creates an impression that the forecast values beyond $t$ is a line.

%\begin{itemize}
%  \item Use \emph{forecast horizons} where necessary
%  \item statsmodels.tsa stands for the statsmodels time series analysis toolbox
%\end{itemize}

The second basic exponential smoothing method is the Holt linear method, which builds on the same path as (\ref{SESformula})--(\ref{SESoptimization}). Holt's linear method is suited for time series involving trend without the presence of seasonality. Hence, this method involves an estimate of the level and linear trend of the time series at a given time point. As a consequence, the Holt linear method involves level and slope parameters $\alpha$ and $\beta$, respectively. These parameters can be optimized using the minimization of the MSE, similarly to what is done in (\ref{SESoptimization}).
Similarly to SES, Holt's linear method is applied by simply calling the function named \texttt{Holt} from \texttt{statsmodels.tsa.api}. In the case where we want to set the parameters $\alpha$ and $\beta$ manually, we can use the options \texttt{smoothing\_level} and \texttt{smoothing\_slope}, respectively.
To improve the forecasting performance of the Holt linear method, the \emph{Holt} function provides an option to select the nature of the trend using the \texttt{exponential} or \texttt{damped} option, as it can be seen in the following excerpt of the Holt forecasting code in Listing \ref{Holt}:

\begin{lstlisting}[language=Python]
fit1 = Holt(series).fit(smoothing_level=0.5, smoothing_slope=0.5, optimized=False)
#------------------------------------
fit2 = Holt(series, exponential=True).fit(smoothing_level=0.2,
            smoothing_slope=0.4, optimized=False)
#------------------------------------
fit3 = Holt(series, damped=True).fit()
\end{lstlisting}
Obviously, the default selection of the trend in the first model (see first line in this excerpt) is the linear trend. More details on the different type of trends and the corresponding mathematical adjustment, see \href{https://www.statsmodels.org/stable/generated/statsmodels.tsa.holtwinters.Holt.html}{https://www.statsmodels.org/stable/generated/statsmodels.tsa.holtwinters.Holt.html}

\begin{figure}[htp]
\subfigure[Data and forecast plots]{\includegraphics[width=4.3cm]{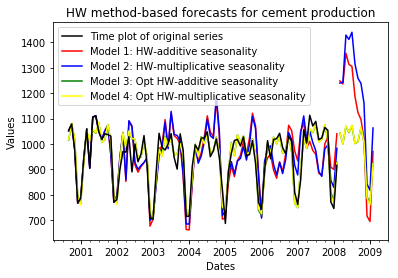}}
\subfigure[Time plots of the residuals]{\includegraphics[width=4.3cm]{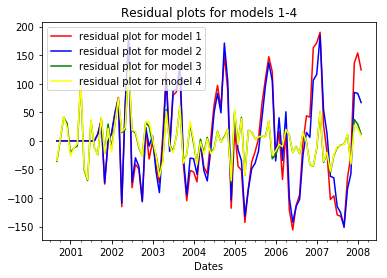}}
\subfigure[Values of parameters and MSE]{\includegraphics[width=6.5cm]{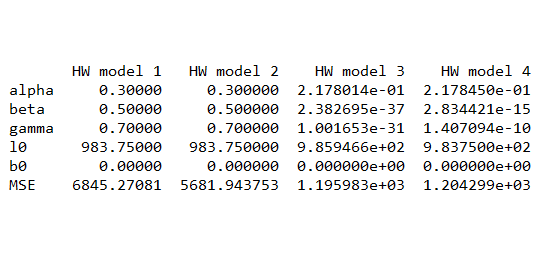}}\\
\subfigure[ACF errors - model 1]{\includegraphics[width=3.71cm]{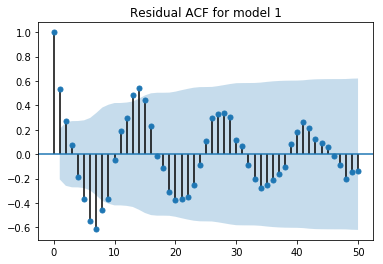}}
\subfigure[ACF errors - model 2]{\includegraphics[width=3.71cm]{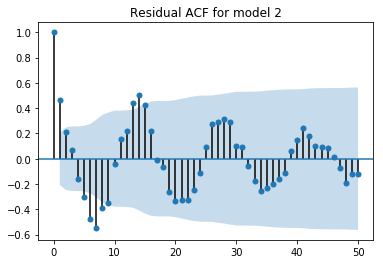}}
\subfigure[ACF errors - model 3]{\includegraphics[width=3.71cm]{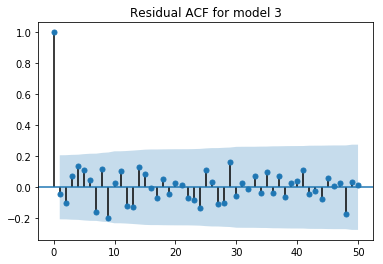}}
\subfigure[ACF errors - model 4]{\includegraphics[width=3.71cm]{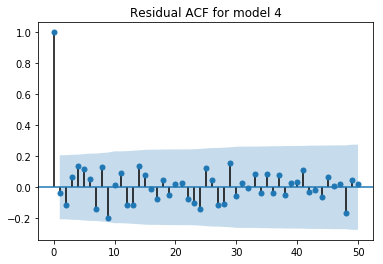}}
\caption{All the results here are generated with the codes in Listing \ref{ExponentialSmoothing}, where (a) is obtained from \texttt{ExponentialSmoothing} and (b) results from errors calculated with values extracted from the results from the \texttt{ExponentialSmoothing} function. The table in (c) is  obtained by applying \texttt{mean\_squared\_error} from \texttt{sklearn.metrics} to forecast values extracted from the results from  \texttt{ExponentialSmoothing}. The latter results can also be obtained straightforwardly from the formulas of the MSE. As for the second row, the plots there are generated with \texttt{plot\_acf} from \texttt{statsmodels}. The data set results form cement production in Australia.}\label{graphHW}
\end{figure}

Finally, we now present the Holt-Winter forecasting method, which is suitable for time series involving both trend and seasonality. Hence, in addition to the level and trend components needed in the Holt linear method (design only for the case where trend in present in our time series), a seasonal components is needed. The seasonal component also comes with its parameter generally denoted by $\gamma$. As it should be the case for the previous two methods, all the parameters are required to be real numbers from the interval $[0, \; 1]$. Since the Holt-Winter method is more general than the SES and LES, the corresponding function from \texttt{statsmodels.tsa.api} is labelled as \texttt{ExponentialSmoothing}.

\begin{lstlisting}[language=Python]
fit1 = ExponentialSmoothing(series, seasonal_periods=12, trend='add',
                            seasonal='add').fit(smoothing_level = 0.3,
                            smoothing_slope=0.5,  smoothing_seasonal=0.7)
\end{lstlisting}
As we can see from this excerpt of the corresponding code in Listing \ref{ExponentialSmoothing}, besides the parameters $\alpha$, $\beta$, and $\gamma$, represented here by \texttt{smoothing\_level}, \texttt{smoothing\_slope}, and \texttt{smoothing\_seasonal},  which can be fixed or optimized as in the previous two exponential smoothing methods, we have the nature of the trend and seasonality, which can be additive or multiplicative. Clearly, the term \emph{add} (resp. \emph{mul}) is used for additive (resp. multiplicative) trend or seasonality. More details on these concepts can be found in \cite[Section 2]{Zemkoho2021}.

We use the code in Listing \ref{ExponentialSmoothing} to generate the results in Figure \ref{graphHW}, which clearly show that the optimized models 3 and 4 are the best, with the 3rd one with additive trend and seasonality being slightly better. The ACF of the residuals from each method are also included in the code, to further evaluate the performance of each method. It is clear that the residuals for models 1 and 2 retain the seasonality present in the original data set. On the other hand Figures \ref{graphHW}(b), (f), and (g) just confirm that residuals seem relatively random.

\section{ARIMA methods}\label{ARIMA methods}

\subsection{Preliminary tools}
As we have seen so far, the ACF plot can play an important role in showing that a time series is seasonal and also in assessing the accuracy of a forecasting method (mainly via the white noise concept). In this section, we are going to see how the ACF can also be helpful in assessing a few other properties relevant to the ARIMA method, namely, in assessing stationarity and the identification of an ARIMA model. However, to strengthen the capacity of the ACF in this role, we now introduce the concept of \emph{partial autocorrelation function} (PACF), which is used  to measure the degree of association between observations at time lags $t$ and $t-k$ (i.e., $Y_t$ and $Y_{t-k}$, respectively) when the effects of other time lags,  $1, \ldots, k-1$, are removed. Hence, partial autocorrelations calculate true correlations between $Y_t$, $Y_{t-1}$, \ldots, $Y_k$ and can therefore be obtained using a regression formula on these terms, while proceeding as in the least square approach in (\ref{SESoptimization}) or the concept of maximum likelihood estimation, which is more common in this case \cite{HyndmanEtAl2018}.

To get a good flavour of how the PACF can be applied, let us use it to further illustrate white noise in combination with ACF.
Similarly to the ACF, as shown in Subsection \ref{Basic data analysis tools}, the PACF can be plotted by simply applying the function \texttt{plot\_pacf} from \texttt{statsmodels.graphics.tsaplots}. The code in Listing \ref{WhiteNoiseModel} generates the AFC and PACF for an example of white noise model. The important thing to note when this code is ran is how the ACF and PACF of a typical white noise model look like; recall that for a model to be statistically while noise, about 95\% of the values of ACF and PACF are within the range $\pm$ $1.96/\sqrt{n}$, where $n$ is the total number of observations. This range is represented by the shadow band that appears in the graphs of both the ACF and PACF.
\begin{figure}[htp]
\hfill
\subfigure{\includegraphics[width=5cm]{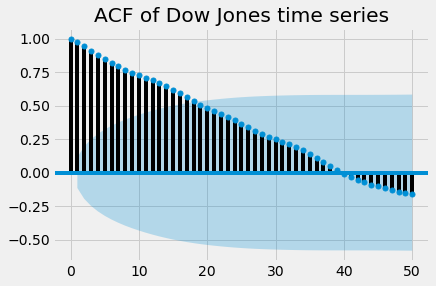}}
\hfill
\subfigure{\includegraphics[width=5cm]{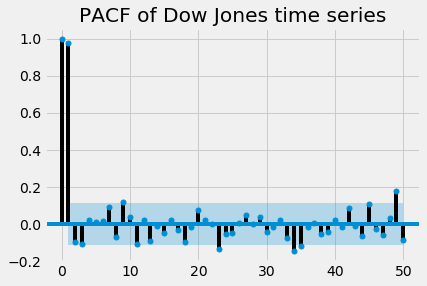}}
\hfill
\caption{Example of non-stationary times series (Dow Jones data from January 1956 to April 1980)}\label{graphNon-Stat}
\end{figure}

We now turn our attention to the concept of \emph{stationarity}, which is at the heart of the development of ARIMA methods. Recall that a time series is  stationary if the distribution of the fluctuations is not time dependent. This is easy to say, but it can be tricky to actually show that a time series is stationary. We try now to provide a few tools that can be helpful in identifying stationarity in a time series. To proceed, we start by stating the following scenarios or specific tools that we are going to rely on to identify whether a time series is stationary or not:
\begin{enumerate}
 \item A white noise time series is stationary;
  \item A time series with trend or seasonality is non-stationary;
  \item A cyclical time series with no trend and no seasonality is stationary;
  \item A non-stationary time series can be detected by means the ACF and PACF;
  \item A unit root test can be used to show that a time series is stationary.
\end{enumerate}

We have just seen how to determine whether a time series is white noise, using the ACF and PACF, which can be plotted with Python using \texttt{plot\_acf} and \texttt{plot\_pacf}, respectively. As for the second item, we already know, see Subsection \ref{Basic data analysis tools}, how to identify trend and seasonality, as well as cyclical patterns, using time plots. There is an interesting way to show that a time series is non-stationary by means of its ACF and PACF plots. Basically, the autocorrelations of a stationary time series drop to zero quite quickly, while those of a non-stationary one can take a significant number of time lags to become zero.  On the other hand, the PACF of a non-stationary time series will typically have a large spike, possibly close to 1, at lag 1. This can clearly be observed in Figure \ref{graphNon-Stat} generated with the code in Listing \ref{AcfPacfPlotDowJones.py}. %, where the left and right pictures correspond to the time series in pictures (a) and (b), respectively. Further details the data and ACF and PACF of the Dow Jones data can be found in the next demo; cf. Listing ...
\begin{figure}[htp]
\subfigure[Original data]{\includegraphics[width=3.4cm]{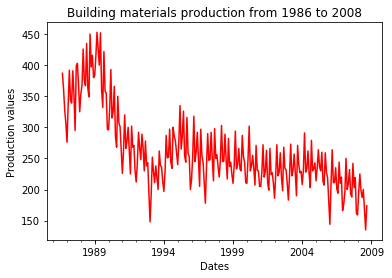}}
\subfigure[ACF plot]{\includegraphics[width=3.4cm]{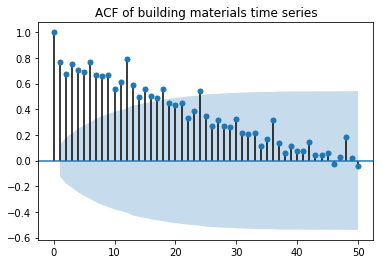}}
\subfigure[PACF plot]{\includegraphics[width=3.4cm]{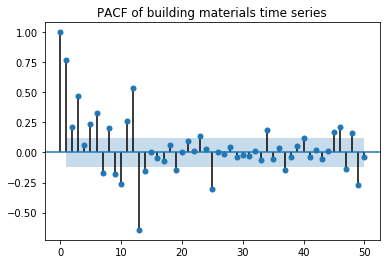}}
\subfigure[ADF test]{\includegraphics[width=4.7cm]{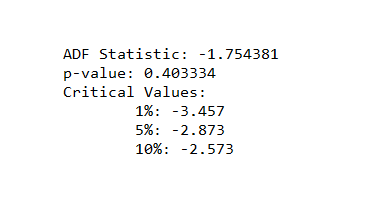}}\\
\subfigure[1st diff]{\includegraphics[width=3.4cm]{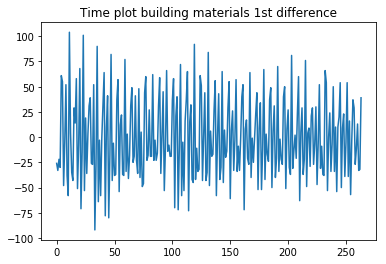}}
\subfigure[1st ACF plot]{\includegraphics[width=3.4cm]{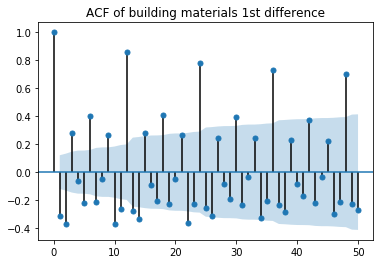}}
\subfigure[1st PACF plot]{\includegraphics[width=3.4cm]{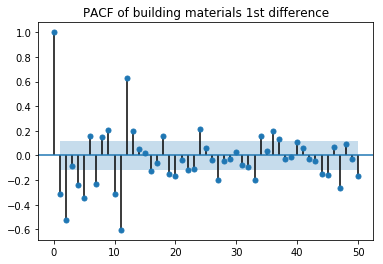}}
\subfigure[1st ADF test]{\includegraphics[width=4.7cm]{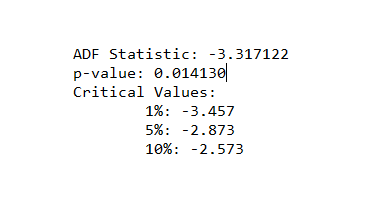}}\\
\subfigure[1st+Sea diff]{\includegraphics[width=3.4cm]{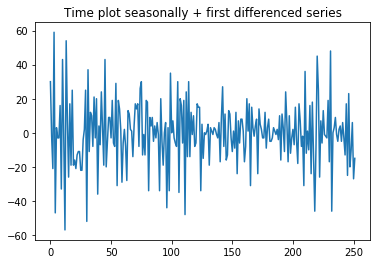}}
\subfigure[1st+Sea ACF plot]{\includegraphics[width=3.4cm]{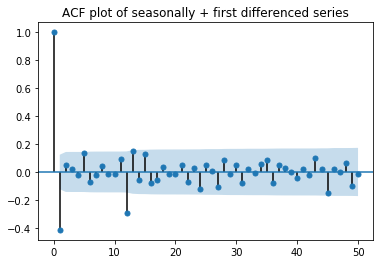}}
\subfigure[1st+Sea PACF plot]{\includegraphics[width=3.4cm]{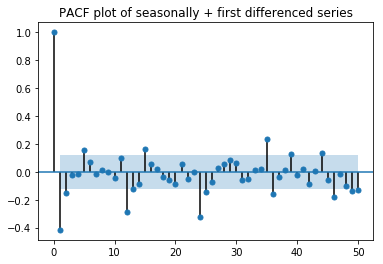}}
\subfigure[1st+Sea ADF test]{\includegraphics[width=4.7cm]{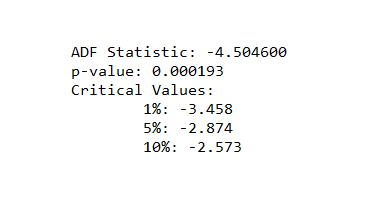}}
\caption{The original time series is obviously not stationary (see first row); after first differening, trend is removed but the new series is still not stationary as it is seasonal (see second row). It is after both first and seasonal differencing that we obtain a stationary time series (see 3rd row).}\label{graphADF}
\end{figure}

Ultimately, if the first four points above cannot help to make a definite decision on the stationarity or non-stationarity of a time series, then we can proceed with a \emph{unit root test}. It is important to say beforehand that this is not a magic solution to demonstrate stationarity, as there are various types of unit root tests, which can sometimes provide contradictory results. The version of the unit root test that we consider here is the Augmented Dickey-Fuller (ADF) test \cite{Dickey-Fuller}, which assesses the null hypothesis that a unit root is present in a time series sample.

A simple understanding of the ADF test that is relevant to us is that it generates a number of statistics that we are going to present next.
To generate these statistics, the function \texttt{adfuller} from \texttt{statsmodels.tsa.stattools} can be applied to our data set. This function simply takes in the values of the time series, as it can be seen in the example used in the code present in Listing \ref{ADFtest.py}, which is used to generate the results in Figure \ref{graphADF} from three different scenarios. Considering some building material production data from Australia,  the first row of Figure \ref{graphADF} presents the time, ACF, and PACF plots, respectively, as well as the statistics generated by the ADF test.

The ADF test (see last column of Figure \ref{graphADF}) generates three key categories of statistics. First, we have the ADF statistics itself, which needs to be negative and subsequently would need to be less than the 1\% critical value to confirm the strength of stationary if additionally, the P-value is at least less than the threshold value of $0.05$. We can clearly see from Figure \ref{graphADF} how the ADF test helps to confirm that we go from a series, where the original and first differenced series are non-stationary to a stationary time series when first and seasonal differencing are done.

\subsection{Models and selection process}
%Three key topics to introduce at the beginning
%\begin{enumerate}
%  \item PACF and ACF?
%  \item White noise (necessary? I think no, as the plan is to have it in the Exp Smth Chapter)
%  \item Stationarity
%\end{enumerate}
%\begin{itemize}
%  \item In practice, ARIMA and Exp Smth are the mostly used methods in forecasting
%  \item Standard stationary time series are: white noise and a time series with cyclic behaviour (but with no trend or seasonality)
%  \item The unit root tests (see Hyndmann et al.) to check whether a time series is stationary
%\end{itemize}
%\begin{enumerate}
%  \item Start to difference while aiming to identifying AR(p) or MA(q) after having ensured stationarity holds
%\end{enumerate}
A simplistic though seemingly nice way to introduce the ARIMA methods is to think about the corresponding forecast function as a polynomial  $f : \mathbb{R}\rightarrow \mathbb{R}$ defined by
\[
f(x):= a_0 + a_1 x + a_2 x^2 + \ldots a_p x^p,
\]
where $p$ is the order of the polynomial and $a_0$, $a_1$, \ldots, $a_p$ are its coefficients. To get a complete description of this polynomial, we need to start by identifying the order $p$, which determines the number of the coefficients $a_0$, $a_1$, \ldots, $a_p$, which can then be subsequently calculated. This is approximately what is done to build an ARIMA model. To make things a bit precise, let us consider a non-seasonal ARIMA$(p, d, q)$ model
\begin{equation}\label{ARIMAgeneral}
(1-\phi_1B - \ldots -\phi_pB^p)(1-B)^dY_t = c+ (1-\theta_1B -\ldots -\theta_qB^q)e_t,
\end{equation}
where $B^kY_t := Y_{t-k}$ corresponds to the backshift notation. Here, the vector $(p, d, q)$ represents the order of the model, and $\phi_i$, $i=1, \ldots, p$ and $\theta_j$, $j=1, \ldots, q$ are parameters/coefficients of the model. Algorithm \ref{algorithm 1} summarizes the building process of an ARIMA model, including the forecasting step.
\begin{algorithm}[htp!]\label{algorithm 1}
\caption{Sketch of the process of identifying and applying an ARIMA model for forecasting}
\label{algorithm 1}
\begin{algorithmic}
\STATE \textbf{Step 1}: Preliminary identification of the order of the model $(p, d, q)$ using ACF and PACF.
\STATE \textbf{Step 2}: Confirmation of the order by means of the Akaike Information Criterion (AIC).
\STATE \textbf{Step 3}: Determine the coefficients $c$, $\phi_1$, \ldots, $\phi_p$, and $\theta_1$, \dots, $\theta_q$.
\STATE \textbf{Step 4}: Apply the order to forecast and assess the quality of the model by means of some statistics.
\end{algorithmic}
\end{algorithm}

\begin{figure}[htp]
\begin{center}
\subfigure{\includegraphics[width=4cm]{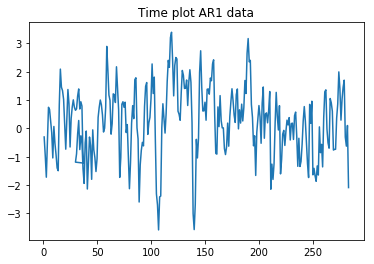}}
\subfigure{\includegraphics[width=4cm]{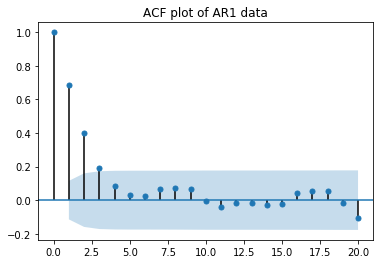}}
\subfigure{\includegraphics[width=4cm]{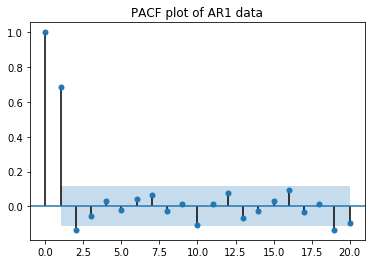}}\\
\subfigure{\includegraphics[width=4cm]{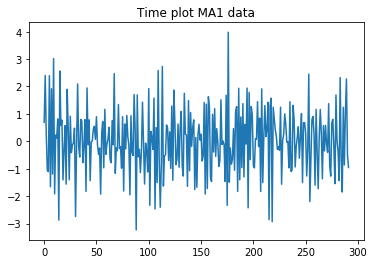}}
\subfigure{\includegraphics[width=4cm]{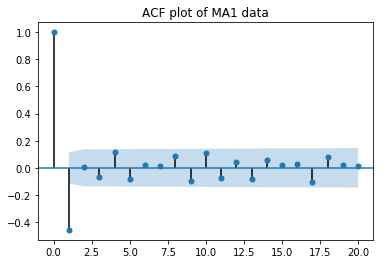}}
\subfigure{\includegraphics[width=4cm]{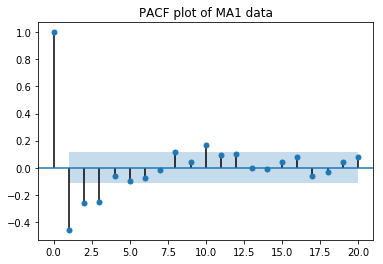}}
\end{center}
\caption{The first row presents the time, ACF, and PACF plots of an artificially generated autoregressive model of order 1. The second row presents analogous graphs for an artificially generated moving average of order 1.}\label{graphAR1&MA1}
\end{figure}

In Step 1, determining $d$ is the most obvious thing, as it simply results from whether we need to do differencing or not to ensure that our time series is stationary. If no differencing is needed, then $d=0$. Otherwise $d\geq 1$ simply represents how many times differencing is needed to obtain stationarity. $p$ and $q$ are much more trickier to obtain. An initial approximation of these numbers can be derived from the ACF and PACF plots. To proceed, let us first present graphs of the ACF and PACF of pure autoregressive and moving average models
\begin{equation}\label{AR&MA}
\mbox{AR}(p) = \mbox{ARIMA}(p, 0, 0) \;\; \mbox{ and } \;\; \mbox{MA}(q) = \mbox{ARIMA}(0, 0, q),
\end{equation}
respectively, for an artificially generated time series example.
AR$(p)$  is obtained if the ACF of this time series is exponentially decaying or sinusoidal and there is a significant spike at lag $p$ in the PACF, but no larger one beyond lag $p$. As for the pure moving average model MA$(q)$, the PACF is exponentially decaying or sinusoidal; and there is a significant spike at lag $q$ in the ACF plot, but not larger than one beyond lag $q$. Of course, in the case of a non-stationary time series, these observations would be made on ACF and PACF of ``sufficiently differenced'' (in the sense of leading to stationarity) data. The graphs in Figure \ref{graphAR1&MA1} show an AR$(1)$ and a MA$(1)$ in the first and second row, as generated by Listings \ref{AcfPacfAR1model.py} and \ref{AcfPacfMA1model.py}, respectively.

\begin{figure}[htp]
\begin{center}
\subfigure{\includegraphics[width=4.85cm]{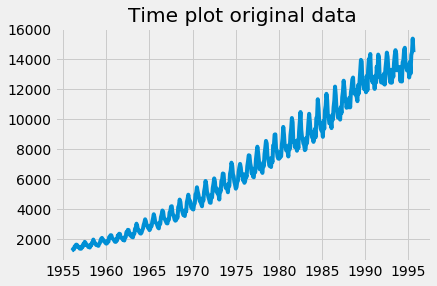}}
\subfigure{\includegraphics[width=4.85cm]{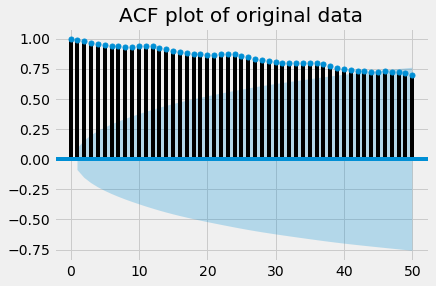}}
\subfigure{\includegraphics[width=4.85cm]{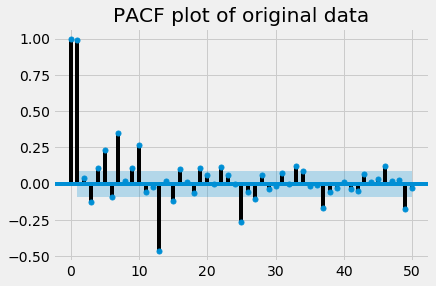}}\\
\subfigure{\includegraphics[width=4.85cm]{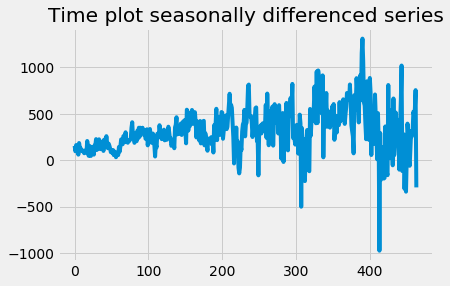}}
\subfigure{\includegraphics[width=4.85cm]{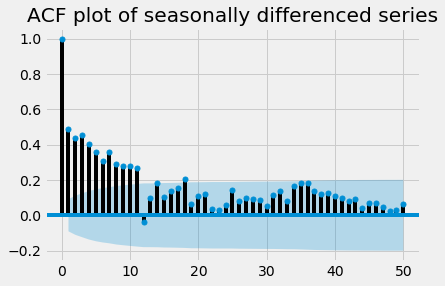}}
\subfigure{\includegraphics[width=4.85cm]{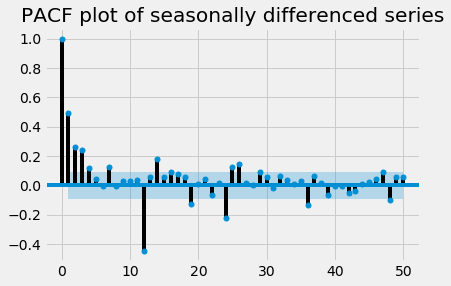}}\\
\subfigure{\includegraphics[width=4.85cm]{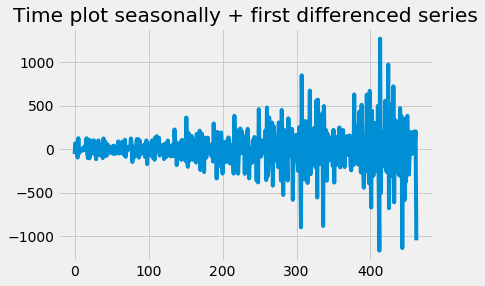}}
\subfigure{\includegraphics[width=4.85cm]{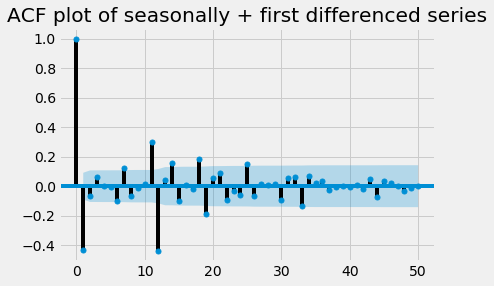}}
\subfigure{\includegraphics[width=4.85cm]{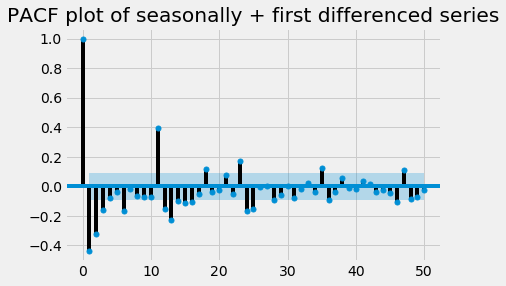}}
\end{center}
\caption{These graphs generated from Listing \ref{FirstSeasonalDifference.py} present the changes in the electricity demand times series data in Figure \ref{graph0}(b), going  from the original data and its ACF and PACF plots (first row), passing by the first difference (second row) to the graphs resulting from first and seasonal differencing (third row).}\label{graphFirstSeasonalDifference}
\end{figure}

Considering the fact that the approach in Step 1 can only enable the estimation of pure AR and MA models, we need a way to check whether our series exhibits a more general ARIMA$(p, d, q)$ model with $p>0$ and $q>0$ simultaneously. The AIC, which is a function of $p$ and $q$ can help us to check whether there is a model better than the one obtained from Step 1.  The smaller the AIC, the better the model is. To proceed, we can use the code in Listing \ref{AutomaticOrderARIMA.py}, which runs through a combinations of values $p$, $d$, and $q$  from the interval $[0, 2]$ to identify the order $(p, d, q)$ with best AIC. %The search interval can be expanded to include values greater than 2.
For the selection of $d$, it is straightforward  to use the process described above, repeating the differencing as necessary to get best statistics from the ADF test based on the code in Listing \ref{ADFtest.py}.
\begin{figure}[htp]
\begin{center}
\subfigure{\includegraphics[width=14cm]{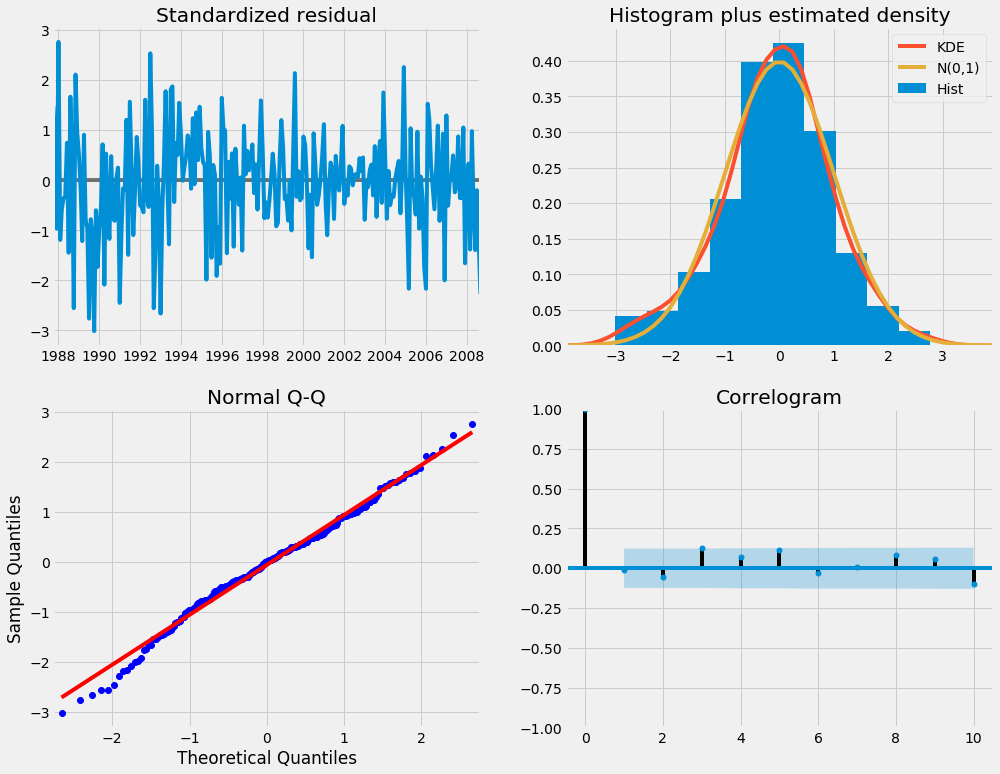}}\\
\subfigure{\includegraphics[width=7.2cm]{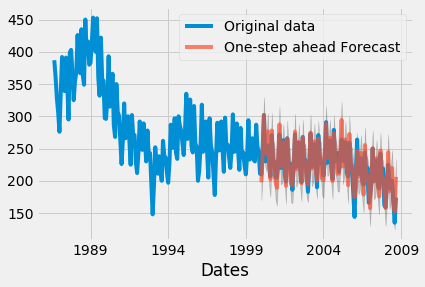}}
\subfigure{\includegraphics[width=6.7cm]{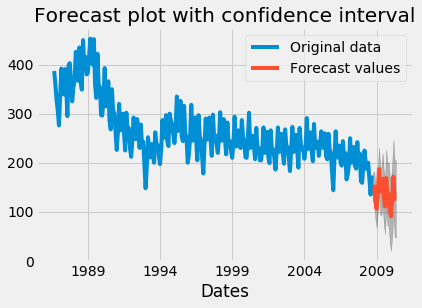}}
\end{center}
\caption{Summary of graphical results obtained by running the SARIMAX$(1, 1, 1)(0, 1, 1)_{12}$ model using the code from Listing \ref{SARIMAX.py} on  building material time series from 1986 to 2008 in Australia. The first four graphs assess the accuracy of the method, with (1) the residual plot, (2) the distribution of the error (close to a normal distribution), (3) the normal Q--Q plot, which compares randomly generated and independent standard normal data on the vertical axis to a standard normal population on the horizontal axis (the closest the data points are to a line suggests that the data are normally distributed), and (4) the correlogram for checking randomness in the residual. The last row shows the one-step forecasts on a section of the data for some visual assessment of accuracy, as well as the out-of-sample future forecasts over a 20 step horizon.}\label{graphArimaResult}
\end{figure}

In terms of the content of the code in Listing \ref{AutomaticOrderARIMA.py}, its main feature is the \texttt{ARIMA} function  from \texttt{statsmodels}. This function is also going to be used for Step 4 of Algorithm \ref{algorithm 1}, but one of its most interesting features is that it also generates other important information such as the AIC of the corresponding model. However, in the context of Listing \ref{AutomaticOrderARIMA.py}, its main role is to print and compare the AIC to identify the best model. When the most suitable values of the order $(p, d, q)$ have  been identified, the \texttt{ARIMA} function can then be applied,  using this order, to generate the forecasts, as it is done for the example in Listing \ref{ARIMA.py}. Running the code generates forecast plots and some important statistics, including the AIC of the model and the corresponding coefficients/parameters $\phi_i$, $i=1, \ldots, p$ and $\theta_j$, $j=1, \ldots, q$ as described in the equation in (\ref{ARIMAgeneral}).

So far, we have considered only time series that are not necessarily seasonal. In the seasonal case, the process is the same, except that the seasonal order $(P, D, Q)$ and periodicity $s$ have to be provided, as indicated in the general model
\begin{equation}\label{SeasonalARIMA}
\mbox{ARIMA}(p, d, q)(P, D, Q)s.
\end{equation}
The first key difference between the non-seasonal (\ref{ARIMAgeneral}) and seasonal (\ref{SeasonalARIMA}) ARIMA models is the parameter $s$, which represents the number of time periods per season shown by the time series; for example, for the seasonal time series examples that we have covered so far (see, e.g., Figures \ref{graph1} and \ref{graphAutoCorrelation}), the patterns repeat themselves every 12 months - hence, $s=12$ in those cases. Similarly to $d$ in (\ref{ARIMAgeneral}), $D$ in (\ref{SeasonalARIMA}) represents the number of seasonal difference needed to remove seasonality in the time series. Furthermore, the pure seasonal autoregressive and moving average models
\begin{equation}\label{SeasonalAR&MA}
 \mbox{ARIMA}(p, 0, 0)(P, 0, 0)_s \;\; \mbox{ and } \;\; \mbox{ARIMA}(0, 0, q)(0, 0, Q)_s,
\end{equation}
respectively, can be obtained in a way similar to  (\ref{AR&MA}) by looking whether the patterns from the ACF and PACF plots approximately repeat themselves after $s$ time lags. For example, using the code in Listing \ref{FirstSeasonalDifference.py}, first order differencing and seasonal differencing can be done to remove trend and seasonality from the electricity data from Figure \ref{graph1}(b). Based on the reference ACF and PACF plots in Figure \ref{graphAR1&MA1} (second row), this leads to the model  \mbox{ARIMA}$(0, 1, 1)(0, 1, 1)_{12}$ in Figure \ref{graphFirstSeasonalDifference} (third row), as a MA(1) pattern approximately repeats itself from the 12th time lag.

After this identification trial for a seasonal model based on the ACF and PACF plots, one can then proceed with the automatic identification process similar to  the one introduced above (see Listing \ref{AutomaticOrderARIMA.py}), while using the corresponding seasonal code available in Listing \ref{AutomaticOrderSARIMAX.py}. Similarly, when the best seasonal order $(p, d, q)(P, D, Q)$ with the corresponding number of time periods per season ($s$) has been identified, the seasonal ARIMA function (\texttt{SARIMAX}) also from \texttt{statsmodels} (see Listing \ref{SARIMAX.py}) can be used to generate the forecasts. Running SARIMAX with the code available in Listing \ref{SARIMAX.py} applied on building material time series from 1986 to 2008 in Australia, we get the graphs in Figure \ref{graphArimaResult} together with a number of statistics assessing the quality of the model and the results.
%Instead of using the function \texttt{ARIMA}, one could apply \texttt{ARIMA}(p, d, q)(0, 0, 0).
%\begin{itemize}
%  \item Highlight the fact that the SARIMAX is better than the ARIMA for the building material used there.
%\end{itemize}

%The remaining important things
%\begin{enumerate}
%  \item How to demonstrate stationarity (\emph{ACF and PACF} can be quite useful for initial guess)--differencing, etc.
%  \begin{itemize}
%    \item The data may follow an ARIMA(p, d, 0) model if the ACF and PACF plots of the differenced data show the following patterns:
%  \end{itemize}
%  \item How to identify the order or the model (two options are possible here):
%  \begin{itemize}
%    \item Using the \emph{ACF and PACF} to make a first guess
%    \item Use the automatic identification tool to check/evaluate
%  \end{itemize}
%  \item When the order of the model is known, all that remains is to use the order to forecast with ARIMA or SARIMAX
%\end{enumerate}

\section{Regression-based forecasting method}\label{Regression-based forecasting method}
The particularity of the method that we are going to discuss here is that it is \emph{explanatory}, in comparison to the previous ones, which are blackbox  methods. A regression model exploits potential relationships between the main (dependent) variable and other (independent) variables. We focus our attention here on the simplest and most commonly used  relationship, which is the \textit{linear regression}:
\begin{equation}\label{linearReg}
 Y= b_o + b_1X_1 + \ldots +b_kX_k +e,
\end{equation}
where $Y$ is the dependent variable, $X_1$, \ldots, $X_k$ the independent variables, and $b_0$, $b_1$, \ldots, $b_k$ the coefficients/parameters, where $b_0$ specifically is often called \emph{intercept}.
It is important to start by recalling that a regression model as (\ref{linearReg}) is not a forecasting method by itself; there is a large number of applications of regression models in statistics and econometrics; see, e.g., \cite{Montgomery2021} for a detailed analysis of regression models and some flavour of a sample of applications.

To apply the regression model (\ref{linearReg}) to develop a forecast for a time series $\{Y_t\}$, we assume that it is influenced by other time series $\{X_{it}\}$ for $i=1, \ldots, n$. To have some flavour of this, we consider the \emph{mutual savings bank case study} from \cite[Chapter 6]{MakridakisEtAl1998}, which arose from a concern in the 90's in the United States with deposits in banks getting smaller than withdrawals. Hence, it became  necessary, in order to anticipate on future decision-making, to develop forecasts of end-of-month (EOM) balance for a few months ahead. Considering the fact that the EOM is influenced by the composite AAA bond rates (just labelled as AAA for short) and the US government 3 to 4 year bond rates (Tto4), cf. Figure \ref{graphCorrelation}, a regression model can be built to forecast EOM while considering AAA and Tto4 as independent variables. For some technical reasons (see \cite{MakridakisEtAl1998}), our $Y$ is the first order difference of EOM (denoted by DEOM), and $X_1$, $X_2$, and $X_3$ the AAA, Tto4, and D3to4 (first order difference of Tto4), respectively. Note that historical time series data sets are available for the variables DEOM, AAA, Tto4, and D3to4, and there are some level of relationship between these variables as it can be seen from the scatter plots and correlation matrix in Figure \ref{graphCorrelation}. However, this is not enough to guarantee that the regression model resulting from this relation would be significant.  The analysis of a regression model starts with the evaluation of its overall significance.

\begin{figure}
\begin{center}
  \includegraphics[width=15cm]{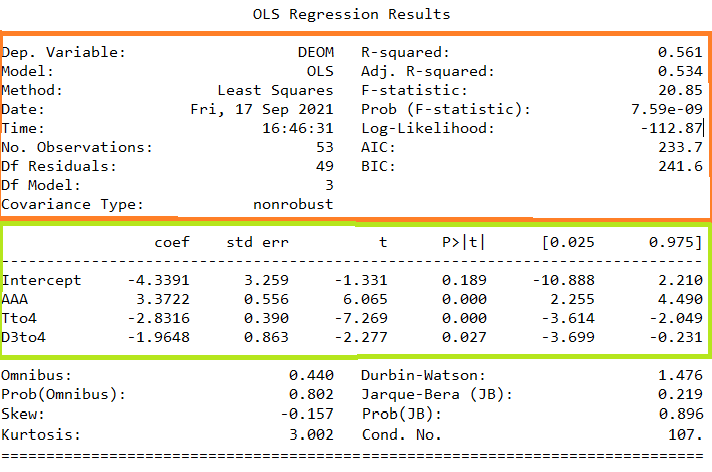}
  \caption{Key statistics to assess the overall and individual significance of a regression model}\label{Regress1}
\end{center}
\end{figure}

For the overall significance of a model, key statistics are the $R^2$ (known as the coefficient of determination) and the $P$-value, which gives the probability of obtaining a $F$ statistic as large as the one calculated for the data set being studied, if in fact the true slope is zero. As the $R^2$ is a number between $0$ and $1$, model (\ref{linearReg}) would be considered to be significant if it is at least greater than $0.50$. Hence, the overall significance of the model increases as $R^2$ grows closer to the upper bound $1$. Furthermore, from the perspective of the $P$-value, a regression model will be said to be significant if the $P$-value is smaller than the conventionally set value of $0.05$; and the significance improves as the $P$-value decreases below this threshold.

Before we expand this discussion further, let us show how the aforementioned statistics can be obtained with Python. Our analysis of a regression model here is based on the \texttt{ols} function from \texttt{statsmodels}, which means \emph{ordinary least squares}, given that the parameters  in (\ref{linearReg}) are computed by the same least square approach introduced for the SES model in (\ref{SESoptimization}). As you can see in the demonstration code in Listing \ref{Regression.py}, it is incredibly easy to use \texttt{ols}. For example, to build the basic model for our above bank case study, what is needed is to start by writing the regression equation
\\[1ex]
\texttt{formula = 'DEOM $\sim$ AAA + Tto4 + D3to4'}, \\[1ex]
recalling that $Y$ = DEOM is the dependent variable and $X_1$ = AAA, $X_2$ = Tto4, and $X_3$ = D3to4, respectively, are the independent variables. The function \texttt{ols} can then be applied as follows to the combination of this formula and the data set to produce the statistics:
\\[1ex]
\texttt{results = ols(formula, data=series).fit()}, \\[1ex]
where \emph{series} corresponds to the container with the time series data sets for the dependent and independent variables. The results from this function (generated with the  code in Listing \ref{Regression.py}) are given in the table in Figure \ref{Regress1}.

%The parameters resulting from these statistics, namely, the coefficients of the regression model (see second statistics table generated by the \texttt{ols} function), which are computed, can then be used to produce the forecasts using the formula presented in the introduction of this Section; cf. {\color{blue}{RegressionForecastsBank.py}}.

\begin{figure}[htp]
\begin{center}
\subfigure{\includegraphics[width=5cm]{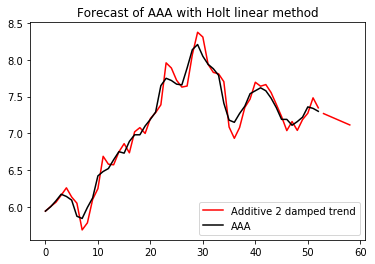}}
\subfigure{\includegraphics[width=5cm]{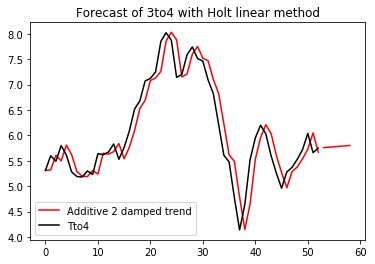}}\\
\subfigure{\includegraphics[width=5cm]{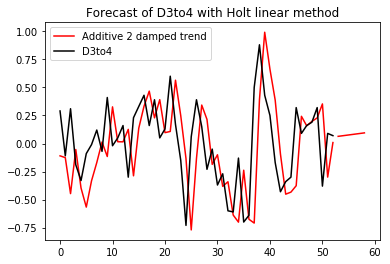}}
\subfigure{\includegraphics[width=5cm]{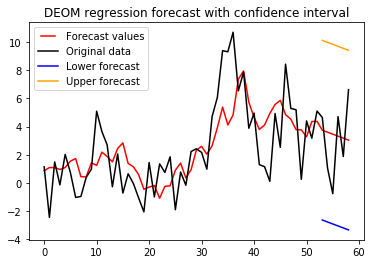}}
\end{center}
\caption{Generating forecasts for the time series involved in this model, i.e., AAA, Tto4, and D3to4 for the independent variables and DEOM for the dependent variable, is quite challenging as none of data sets exhibits a clear pattern. Hence, from the from the exponential smoothing methods covered in Section \ref{Exponential smoothing methods0}, only Holt's linear method is the most suitable, as it enables the calculation of out-of-sample forecasts over a number of time points ahead. An ARIMA method could also be used to generates forecasts for AAA, Tto4, and D3to4.}\label{graphRegRes}
\end{figure}

The orange box in  Figure \ref{Regress1} contains the statistics of overall significance of the model from the corresponding example, where $Y$ is represented by time series DEOM, while $X_1$, $X_2$, $X_3$, and $X_4$, are represented by AAA, Tto4, and D3to4, respectively. It can clearly be seen that the overall model in this example is significant, with an $R^2$ of $0.56$ and a P-value of $7.59\times 10^{-9}$. But the $R^2$ suggests that the significance is not that strong, although the P-value is relatively good from the perspective of the threshold value of $0.05$.

For the individual significance of each variable involved in (\ref{linearReg}), the main statistics is the $P$-value. Similarly to the \emph{F}-test, the \texttt{statsmodels} function  \texttt{ols} (ordinary least squares) generates $P$-values of each $t$-statistic. Each of these $P$-values is the probability of obtaining an absolute value of the $t$-statistics, of a given independent variable, as large as the one calculated for the data, if  the parameter is equal to zero. So, if a $P$-value is small, then the estimated parameter is significantly different from zero. As with
\emph{F}-tests, it is common to conclude that an estimated parameter is significantly different from zero (i.e., \emph{significant}) if the $P$-value is smaller than $0.05$. The 5th column of the green box in Figure \ref{Regress1} gives the $P$-value of each of the 3 independent variables in the example introduced above. Clearly, the significance of AAA, Tto4, and D3to4 is relatively good, as it is less than the threshold value of $0.05$, although that of the latter variable is weaker.

Interestingly, the green box in the table in Figure \ref{Regress1} also provides the coefficients of this example (cf. second column). After we have seen how the function \texttt{ols} can help to generate the key statistics to assess the overall and individual significance of the model, it remains to see how the forecast can actually be derived. To be able to do this, we need the forecasts
\[
G_i = (G_{i1}, \ldots, G_{ik})\; \mbox{ of }\; X_i=(X_{i1}, \ldots, X_{ik}) \;\mbox{ for } \; i = t+ 1, \ldots, t+m.
\]
We can then use each of these forecasts of the independent variables in the expected value that determined the regression-based forecast for the independent variable $Y$ using equation (\ref{linearReg}):
\begin{equation}\label{RegForecast}
 F_i = \hat{Y}_i = G_i \hat{b}\; \mbox{ for } \; i = t+ 1, \ldots, t+m,
\end{equation}
where the forecasts $G_i$ of each independent variable can be obtained by any method that is most suitable. Applying (\ref{RegForecast})to our example above (see Listing \ref{RegressionForecast.py}), we obtain the results in Figure \ref{graphRegRes}.

To conclude this section, some quick comments are in order. First, one of the typical preliminary step when building a regression model is to conduct a  correlation analysis (e.g., scatter plots, correlation matrix),  which can be done using tools that we have discussed in Subsection \ref{Basic data analysis tools}. This can be done here with matrix scatter plots and correlation tables; see Figure \ref{graphCorrelation}.
Also, often to improve an initial model as in (\ref{linearReg}) or resulting forecasting accuracy (\ref{RegForecast}), a careful selection process of variables or features of the data sets can be done. Finally, the term prediction is usually confused with that of forecast. \emph{Prediction} is much more broad, as it includes tasks such as predicting the result of a soccer game or an election, where only characteristics of the player of each team (soccer) or surveys from voters (election) not necessarily historical data can be used. Further details on these topics can be found in \cite{HyndmanEtAl2018,MakridakisEtAl1998,Zemkoho2021} and references therein.

\section{Conclusion}
This paper puts together a set of Python--based mostly off--the--shelf tools to develop forecasts for time series data using basic statistical forecasting methods, namely, exponential smoothing, ARIMA, and regression methods. It is important to mention that for each forecasting method and analysis tool described in this paper,  there could be multiple Python approaches available, to undertake them, across different Python--based platforms. Secondly, within many packages, there could also be various ways to do the same thing. So, when using the material presented here, it will be useful to have a look at the most recent updates on the corresponding packages' websites (see corresponding links provided in Section \ref{Preliminaries on Python and data analysis}) for other possible ways to conduct specific analysis or for the most recent updates on possible improvements to these tools.

\section*{Acknowledgement}
The lecture notes \cite{Zemkoho2021} (based on the textbooks \cite{HyndmanEtAl2018,MakridakisEtAl1998}), which have served as base for the mathematical background of the data analysis and forecasting tools discussed in  this paper, have been developed and refined over the years thanks to contributions from many colleagues from the Southampton OR Group, in particular, I would like to mention Russell Cheng and Honora Smith for preparing and delivering the Forecasting course for many years, until the 2013-14 academic year.

The author would like to thank the referee and the guest editor for their constructive feedback, which led to improvements in the presentation of the paper.

This work is supported by the EPSRC grant with reference EP/V049038/1 and the Alan Turing Institute under the
EPSRC grant EP/N510129/1.

 %It is also important to recall the textbooks for \cite{Zemkoho2021} are \cite{HyndmanEtAl2018,MakridakisEtAl1998}.
%
%. The framework for putting together the Python framework for the delivery of this course was possible in large part thanks to their effort. % largely  The author is grateful to Russell Cheng and Honora Smith  for developing initial Excel-based versions of the lecture notes \cite{Zemkoho2021}, which have serve as based of the mathematical content of the material presented that Also add an acknowledgement for Russel Cheng who taught the ....
%The material The author is also indebted to any colleagues in the Southampton OR Group who have provided any feedback over the years that help in improving the content of the Forecating course content and its delivery.

\section*{Conflict of interest statement}
To the best of his knowledge, the author is not aware of any  conflict of interest.

\section*{Data availability statement} All the data sets used for the illustrations in this paper are based on the book \cite{MakridakisEtAl1998}; all the data sets related to this book are available online: \href{https://cloud.r-project.org/web/packages/fma/index.html}{cloud.r-project.org/web/packages/fma/index.html}. As for the specific times series from this data base used in this paper, they are available via the following link, together with all the \texttt{py} files associated to the codes in the appendix:  \href{https://github.com/abzemkoho/forecasting}{github.com/abzemkoho/forecasting}.
%All the data sets used for the illustrations in this paper can be accessed via this link:\\

\newpage
\appendix

\counterwithin{lstlisting}{section}
\section{Codes on preliminaries on Python and data analysis}

\begin{lstlisting}[language=Python, caption={TimePlot}, label={TimePlot}]
from pandas import read_excel
from matplotlib import pyplot
series = read_excel('ClayBricks.xls', sheetname='BRICKSQ', header=0,
              index_col=0, parse_dates=True, squeeze=True)
series.plot()
pyplot.show()
\end{lstlisting}

\begin{lstlisting}[language=Python, caption={SeasonalPlot}, label={SeasonalPlot}]
from pandas import read_excel
import matplotlib.pyplot as plt
import numpy as np
series = read_excel('ClayBricks.xls',
              sheetname='SeasonalData', header=0, index_col=0, parse_dates=True, squeeze=True)
x = np.array([0,1,2,3,4,5,6,7,8,9,10,11])
months = ['Jan','Feb','Mar','Apr', 'May', 'Jun', 'Jul', 'Aug', 'Sep', 'Oct', 'Nov', 'Dec']
plt.xticks(x, months)
plt.plot(x, series)
plt.show()
\end{lstlisting}

\begin{lstlisting}[language=Python, caption={LogTransform}, label={LogTransform}]
from pandas import read_excel
from pandas import DataFrame
from numpy import log
from matplotlib import pyplot
series = read_excel('Electricity.xls',
              sheetname='Data', header=0,
              index_col=0, parse_dates=True, squeeze=True)
dataframe = DataFrame(series.values)
dataframe.columns = ['electricity']
dataframe['electricity'] = log(dataframe['electricity'])
pyplot.figure(1)
# line plot
pyplot.subplot(211)
pyplot.plot(dataframe['electricity'])
# histogram
pyplot.subplot(212)
pyplot.hist(dataframe['electricity'])
pyplot.show()
\end{lstlisting}

\begin{lstlisting}[language=Python, caption={SqrtTransform}, label={SqrtTransform}]
from pandas import read_excel
from pandas import DataFrame
from numpy import sqrt
from matplotlib import pyplot
series = read_excel('Electricity.xls',
              sheetname='Data', header=0,
              index_col=0, parse_dates=True, squeeze=True)
dataframe = DataFrame(series.values)
dataframe.columns = ['electricity']
dataframe['electricity'] = sqrt(dataframe['electricity'])
pyplot.figure(1)
# line plot
pyplot.subplot(211)
pyplot.plot(dataframe['electricity'])
# histogram
pyplot.subplot(212)
pyplot.hist(dataframe['electricity'])
pyplot.show()
\end{lstlisting}

\begin{lstlisting}[language=Python, caption={CalendarAdjustment}, label={CalendarAdjustment}]
from pandas import read_excel
from matplotlib import pyplot
series = read_excel('MilkProduction.xls',
              sheetname='AdjustedData', header=0,
              index_col=0, parse_dates=True, squeeze=True)  # you can include various other parameters
series.plot()
pyplot.show()
\end{lstlisting}

%\begin{lstlisting}[language=Python, caption=Demo 1.5 (Calculating/plotting/representing the autocorrelation), label={autocorrelation1}]
%# create an autocorrelation plot
%from pandas import read_excel
%from matplotlib import pyplot
%from pandas.plotting import autocorrelation_plot
%series = read_excel('Beer.xls', sheet_name='Data', header=0, index_col=0, parse_dates=True, squeeze=True)
%autocorrelation_plot(series)
%pyplot.show()
%\end{lstlisting}
%
%
%\begin{lstlisting}[language=Python, caption=Demo 1.6 (Trend estimation/calculation/representation)]
%from pandas import read_excel
%from matplotlib import pyplot
%series = read_excel('HouseSales.xls',
%              sheetname='MAData', header=0,
%              index_col=0, parse_dates=True, squeeze=True)  # you can include various other parameters
%series.plot()
%pyplot.show()
%\end{lstlisting}
%
\begin{lstlisting}[language=Python, caption={DecompositionAdditive}, label={DecompositionAdditive}]
from pandas import read_excel
from matplotlib import pyplot
from statsmodels.tsa.seasonal import seasonal_decompose
series = read_excel('HouseSales.xls', sheetname='Data', header=0,
              index_col=0, parse_dates=True, squeeze=True)
result = seasonal_decompose(series, model='additive')
#result = seasonal_decompose(series, model='multiplicative')
result.plot()
pyplot.show()
\end{lstlisting}

\begin{lstlisting}[language=Python, caption={DecompositionMultiplicative}, label={DecompositionMultiplicative}]
from pandas import read_excel
from matplotlib import pyplot
from statsmodels.tsa.seasonal import seasonal_decompose
series = read_excel('AirlineSales.xls', sheetname='Data', header=0,
              index_col=0, parse_dates=True, squeeze=True)
#result = seasonal_decompose(series, model='additive')
result = seasonal_decompose(series, model='multiplicative')
result.plot()
pyplot.show()
\end{lstlisting}

\begin{lstlisting}[language=Python,  caption={ScatterPlot}, label={correlation}]
import pandas as pd
import numpy as np
import matplotlib.pyplot as plt
from pandas import read_excel

#Reading the data from Excel files
series1 = read_excel('JapaneseCars.xls', sheet_name='Data', usecols=[0], header=0,
                     squeeze=True, dtype=float)
series2 = read_excel('JapaneseCars.xls', sheet_name='Data', usecols=[1], header=0,
                     squeeze=True, dtype=float)

#Generating the scatter plot for variables
Japanese = {'Mileage': series1, 'Price': series2}
df = pd.DataFrame(Japanese, columns=['Mileage', 'Price'])
plt.scatter(df.Mileage, df.Price)
plt.xlabel('Mileage (MPG)')
plt.ylabel('Price (US$)')
plt.title('Price/Mileage Relationship for Japanese cars')
plt.show()

#Calculating & printing the correlation between the 2 variables
correlval=np.corrcoef(series1, series2)
correlval=correlval[1,0]
print('\nThe correlation between the 2 variables is:{}'.format(correlval))
\end{lstlisting}

\begin{lstlisting}[language=Python, caption={CorrelationMatrix}, label={CorrelationMatrix}]
from pandas import read_excel
import matplotlib.pyplot as plt
import pandas as pd
series = read_excel('Bank.xls', sheet_name='Data3', header=0,
                     squeeze=True, dtype=float)

#Plotting the scatter plots of each variable against the other one
pd.plotting.scatter_matrix(series, figsize=(8, 8))
plt.show()

# Correlation matrix for all the variables, 2 by 2
CorrelationMatrix = series.corr()
print(CorrelationMatrix)
# As in the case of Demo 1.4, corrcoef can be used for the variables in couples.
\end{lstlisting}

\begin{lstlisting}[language=Python, caption={Autocorrelation}, label={Autocorrelation}]
from pandas import read_excel
from matplotlib import pyplot
from pandas.plotting import autocorrelation_plot
from statsmodels.graphics.tsaplots import plot_acf
series1 = read_excel('CementProduction.xls', sheetname='Data', header=0,
              index_col=0, parse_dates=True, squeeze=True)
series2 = read_excel('CementProduction.xls', sheetname='SeasonalData', header=0,
                    index_col=0, parse_dates=True, squeeze=True)
series2.plot(title='Seasonal plots building materials time series')
pyplot.show()

plot_acf(series1, title='ACF plot of building materials time series', lags=60)
pyplot.show()

autocorrelation_plot(series1)
pyplot.show()
\end{lstlisting}

%%%%%%%%%%%%%%%%%%%%%%%%%%%%%%%%%%%%%%%%%%%%%%%%%%%%%%%%%%%%%%%%%%%%%%%%%%%%%%%%%%
%%%%%%%%%%%%%%%%%%%%%%%%%%%%%%%%%%%%%%%%%%%%%%%%%%%%%%%%%%%%%%%%%%%%%%%%%%%%%%%%%%
%%%%%%%%%%%%%%%%%%%%%%%%%%%%%%%%%%%%%%%%%%%%%%%%%%%%%%%%%%%%%%%%%%%%%%%%%%%%%%%%%%
%%%%%%%%%%%%%%%%%%%%%%%%%%%%%%%%%%%%%%%%%%%%%%%%%%%%%%%%%%%%%%%%%%%%%%%%%%%%%%%%%%
%%%%%%%%%%%%%%%%%%%%%%%%%%%%%%%%%%%%%%%%%%%%%%%%%%%%%%%%%%%%%%%%%%%%%%%%%%%%%%%%%%
%%%%%%%%%%%%%%%%%%%%%%%%%%%%%%%%%%%%%%%%%%%%%%%%%%%%%%%%%%%%%%%%%%%%%%%%%%%%%%%%%%
%%%%%%%%%%%%%%%%%%%%%%%%%%%%%%%%%%%%%%%%%%%%%%%%%%%%%%%%%%%%%%%%%%%%%%%%%%%%%%%%%%
%%%%%%%%%%%%%%%%%%%%%%%%%%%%%%%%%%%%%%%%%%%%%%%%%%%%%%%%%%%%%%%%%%%%%%%%%%%%%%%%%%
%%%%%%%%%%%%%%%%%%%%%%%%%%%%%%%%%%%%%%%%%%%%%%%%%%%%%%%%%%%%%%%%%%%%%%%%%%%%%%%%%%
%%%%%%%%%%%%%%%%%%%%%%%%%%%%%%%%%%%%%%%%%%%%%%%%%%%%%%%%%%%%%%%%%%%%%%%%%%%%%%%%%%
\section{Codes on exponential smoothing methods}

\begin{lstlisting}[language=Python, caption={ErrorMeasures}, label={ErrorMeasures}]
from pandas import read_excel
from matplotlib import pyplot
AustralianBeer  = read_excel('BeerErrorsData.xls', sheet_name='NF1NF2', usecols = [1],
                             header=0, squeeze=True, dtype=float)
NaiveF1  = read_excel('BeerErrorsData.xls', sheet_name='NF1NF2', usecols = [2],
                      header=0, squeeze=True, dtype=float)
NaiveF2 = read_excel('BeerErrorsData.xls', sheet_name='NF1NF2', usecols=[3],
                     header=0, squeeze=True, dtype=float)

# Joint plot of original data and NF1 forecasts
AustralianBeer.plot(legend=True)
NaiveF1.plot(legend=True)
pyplot.show()

# Joint plot of original data and NF2 forecasts
AustralianBeer.plot(legend=True)
NaiveF2.plot(legend=True)
pyplot.show()

# Evaluating the errors from both NF1 and NF2 methods
Error1 = AustralianBeer - NaiveF1
Error2 = AustralianBeer - NaiveF2
ME1 = sum(Error1)* 1.0/len(NaiveF1)
ME2 = sum(Error2)* 1.0/len(NaiveF2)
MAE1=sum(abs(Error1))*1.0/len(NaiveF1)
MAE2=sum(abs(Error2))*1.0/len(NaiveF2)
MSE1=sum(Error1**2)*1.0/len(NaiveF1)
MSE2=sum(Error2**2)*1.0/len(NaiveF2)

PercentageError1=(Error1/AustralianBeer)*100
PercentageError2=(Error2/AustralianBeer)*100
MPE1 = sum(PercentageError1)* 1.0/len(NaiveF1)
MPE2 = sum(PercentageError2)* 1.0/len(NaiveF2)
MAE1=sum(abs(PercentageError1))*1.0/len(NaiveF1)
MAE2=sum(abs(PercentageError2))*1.0/len(NaiveF2)


#Printing a summary of the errors in a tabular form
print('Summary of errors resulting from NF1 & NF2:')
import pandas as pd
cars = {'Errors': ['ME','MAE','MSE','MPE', 'MAPE'],
        'NF1': [ME1, MAE1, MSE1, MPE1, MAE1],
        'NF2': [ME2, MAE2, MSE2, MPE2, MAE2]
        }
AllErrors = pd.DataFrame(cars, columns = ['Errors', 'NF1', 'NF2'])
print(AllErrors)
\end{lstlisting}

\begin{lstlisting}[language=Python, caption={ACFErrors}, label={ACFErrors}]
from pandas import read_excel
from matplotlib import pyplot
AustralianBeer  = read_excel('BeerErrorsData.xls', sheet_name='NF1NF2', usecols = [1],
                             header=0, squeeze=True, dtype=float)
NaiveF1  = read_excel('BeerErrorsData.xls', sheet_name='NF1NF2', usecols = [2],
                      header=0, squeeze=True, dtype=float)
NaiveF2 = read_excel('BeerErrorsData.xls', sheet_name='NF1NF2', usecols=[3],
                     header=0, squeeze=True, dtype=float)

# Plot for the original data set
AustralianBeer.plot(label='Original data', legend=True)
pyplot.show()

# Evaluating the errors from both NF1 and NF2 methods
Error1 = AustralianBeer - NaiveF1
Error2 = AustralianBeer - NaiveF2

# Plot of the error time series
Error1.plot(label='NF1 error plot', legend=True)
Error2.plot(label='NF2 error plot', legend=True)
pyplot.show()

# Creating an autocorrelation plot
from pandas.plotting import autocorrelation_plot
autocorrelation_plot(Error1)
autocorrelation_plot(Error2)
pyplot.show()
\end{lstlisting}

\begin{lstlisting}[language=Python, caption={ConfidenceInterval}, label={ConfidenceInterval}]
from pandas import read_excel
from matplotlib import pyplot
from numpy import sqrt
AustralianBeer  = read_excel('BeerErrorsData.xls', sheet_name='NF1NF2', usecols = [1],
                             header=0, squeeze=True, dtype=float)
NaiveF1  = read_excel('BeerErrorsData.xls', sheet_name='NF1NF2', usecols = [2],
                      header=0, squeeze=True, dtype=float)
NaiveF2 = read_excel('BeerErrorsData.xls', sheet_name='NF1NF2', usecols=[3],
                     header=0, squeeze=True, dtype=float)



# Evaluating the errors from both NF1 and NF2 methods
Error1 = AustralianBeer - NaiveF1
Error2 = AustralianBeer - NaiveF2
MSE1=sum(Error1**2)*1.0/len(NaiveF1)
MSE2=sum(Error2**2)*1.0/len(NaiveF2)

LowerForecast1 = NaiveF1 - 1.645*sqrt(MSE1)
UpperForecast1 = NaiveF1 + 1.645*sqrt(MSE1)

LowerForecast2 = NaiveF2 - 1.645*sqrt(MSE2)
UpperForecast2 = NaiveF2 + 1.645*sqrt(MSE2)

# Joint plot of original data and NF1 forecasts
AustralianBeer.plot(label='Original data')
NaiveF1.plot(label='NF1 forecast')
LowerForecast1.plot(label='NF1 lower bound')
UpperForecast1.plot(label='NF1 upper bound')
pyplot.legend()
pyplot.show()

# Joint plot of original data and NF2 forecasts
AustralianBeer.plot(label='Original data')
NaiveF2.plot(label='NF2 forecast')
LowerForecast2.plot(label='NF2 lower bound')
UpperForecast2.plot(label='NF2 upper bound')
pyplot.legend()
pyplot.show()
\end{lstlisting}

\begin{lstlisting}[language=Python, caption={SimpleExpSmoothing}, label={SimpleExpSmoothing}]
from pandas import read_excel
from statsmodels.tsa.api import SimpleExpSmoothing
from matplotlib import pyplot
series = read_excel('ShampooSales.xls', sheetname='Data', header=0,
              index_col=0, parse_dates=True, squeeze=True)

# Simple Exponential Smoothing #

## SES model 1: alpha = 0.1
fit1 = SimpleExpSmoothing(series).fit(smoothing_level=0.1,optimized=False)
fcast1 = fit1.forecast(10).rename(r'$\alpha=0.1$')
# Plot of fitted values and forecast of next 10 values, respectively
fit1.fittedvalues.plot(color='blue')
fcast1.plot(color='blue', legend=True)

## SES model 2: alpha = 0.7
fit2 = SimpleExpSmoothing(series).fit(smoothing_level=0.7,optimized=False)
fcast2 = fit2.forecast(10).rename(r'$\alpha=0.7$')
# Plot of fitted values and forecast of next 10 values, respectively
fcast2.plot(color='red', legend=True)
fit2.fittedvalues.plot(color='red')

## SES model 3: alpha automatically selected by the built-in optimization software
fit3 = SimpleExpSmoothing(series).fit()
fcast3 = fit3.forecast(10).rename(r'$\alpha=%s$'%fit3.model.params['smoothing_level'])
# Plot of fitted values and forecast of next 10 values, respectively
fcast3.plot(color='green', legend=True)
fit3.fittedvalues.plot(color='green')

# Plotting the original data together with the 3 forecast plots
series.plot(color='black', legend=True)
pyplot.show()

#Evaluating the errors
from sklearn.metrics import mean_squared_error
MSE1=mean_squared_error(fit1.fittedvalues, series)
MSE2=mean_squared_error(fit2.fittedvalues, series)
MSE3=mean_squared_error(fit3.fittedvalues, series)

print('Summary of errors resulting from SES models 1, 2 & 3:')
import pandas as pd
cars = {'Model': ['MSE'],
        'SES model 1': [MSE1],
        'SES model 2': [MSE2],
        'SES model 3': [MSE3]
        }
AllErrors = pd.DataFrame(cars, columns = ['Model', 'SES model 1', 'SES model 2', 'SES model 3'])
print(AllErrors)
\end{lstlisting}

\begin{lstlisting}[language=Python, caption={Holt}, label={Holt}]
from pandas import read_excel
from statsmodels.tsa.api import Holt
from matplotlib import pyplot
series = read_excel('ShampooSales.xls', sheetname='Data', header=0,
              index_col=0, squeeze=True)

# Holt model 1: alpha = 0.8, beta=0.2
fit1 = Holt(series).fit(smoothing_level=0.8, smoothing_slope=0.2, optimized=False)
fcast1 = fit1.forecast(12).rename("Holt's linear trend")

fit2 = Holt(series, exponential=True).fit(smoothing_level=0.8, smoothing_slope=0.2, optimized=False)
fcast2 = fit2.forecast(12).rename("Exponential trend")

fit3 = Holt(series, damped=True).fit(smoothing_level=0.8, smoothing_slope=0.2)
fcast3 = fit3.forecast(12).rename("Additive damped trend")

fit4 = Holt(series).fit(optimized=True)
fcast4 = fit4.forecast(12).rename("Additive 2 damped trend")


fit1.fittedvalues.plot(color='blue')
fcast1.plot(color='blue', legend=True)

fit2.fittedvalues.plot(color='red')
fcast2.plot(color='red', legend=True)

fit3.fittedvalues.plot(color='green')
fcast3.plot(color='green', legend=True)

fcast4.plot(color='yellow', legend=True)

series.plot(color='black', legend=True)
pyplot.show()

#Evaluating the errors
from sklearn.metrics import mean_squared_error
MSE1=mean_squared_error(fit1.fittedvalues, series)
MSE2=mean_squared_error(fit2.fittedvalues, series)
MSE3=mean_squared_error(fit3.fittedvalues, series)

print('Summary of errors resulting from SES models 1, 2 & 3:')
import pandas as pd
cars = {'Model': ['MSE'],
        'LES model 1': [MSE1],
        'LES model 2': [MSE2],
        'LES model 3': [MSE3]
        }
AllErrors = pd.DataFrame(cars, columns = ['Model', 'LES model 1', 'LES model 2', 'LES model 3'])
print(AllErrors)
\end{lstlisting}

\begin{lstlisting}[language=Python, caption={ExponentialSmoothing}, label={ExponentialSmoothing}]
from pandas import read_excel
import pandas as pd
from statsmodels.tsa.api import ExponentialSmoothing
from matplotlib import pyplot
series = read_excel('CementProduction.xls', sheetname='Data', header=0,
              index_col=0, parse_dates=True, squeeze=True)

# ===================================
# Holt-Winter method in different scenarios #
# ===================================
# ===================================
# Model 1: Holt-Winter method with additive trend and seasonality
# Here, alpha = 0.3, beta=0.5, gamma=0.7
# ===================================
fit1 = ExponentialSmoothing(series, seasonal_periods=12, trend='add', seasonal='add').fit(smoothing_level = 0.3, smoothing_slope=0.5,  smoothing_seasonal=0.7)
fit1.fittedvalues.plot(color='red')

# ===================================
# Model 2: Holt-Winter method with additive trend and multiplicative seasonality
# Here, alpha = 0.3, beta=0.5, gamma=0.7
# ===================================
fit2 = ExponentialSmoothing(series, seasonal_periods=12, trend='add', seasonal='mul').fit(smoothing_level = 0.3, smoothing_slope=0.5,  smoothing_seasonal=0.7)
fit2.fittedvalues.plot(color='blue')

# ===================================
# Model 3: Holt-Winter method with additive trend and seasonality
# Here, the parameters alpha, beta, and gamma are optimized
# ===================================
fit3 = ExponentialSmoothing(series, seasonal_periods=12, trend='add', seasonal='add').fit()
fit3.fittedvalues.plot(color='green')

# ===================================
# Model 4: Holt-Winter method with additive trend and multiplicative seasonality
# Here, the parameters alpha, beta, and gamma are optimized
# ===================================
fit4 = ExponentialSmoothing(series, seasonal_periods=12, trend='add', seasonal='mul').fit()
fit4.fittedvalues.plot(color='yellow')

print("Forecasting Cement Production with Holt-Winters method")
#=====================================
# Time and forecast plots
#=====================================
series.rename('Time plot of original series').plot(color='black', legend=True)
fit1.forecast(12).rename('Model 1: HW-additive seasonality').plot(color='red', legend=True)
fit2.forecast(12).rename('Model 2: HW-multiplicative seasonality').plot(color='blue', legend=True)
fit3.forecast(12).rename('Model 3: Opt HW-additive seasonality').plot(color='green', legend=True)
fit4.forecast(12).rename('Model 4: Opt HW-multiplicative seasonality').plot(color='yellow', legend=True)
pyplot.xlabel('Dates')
pyplot.ylabel('Values')
pyplot.title('HW method-based forecasts for cement production')
pyplot.show()


#====================================
# Evaluating the errors
#====================================
from sklearn.metrics import mean_squared_error
MSE1=mean_squared_error(fit1.fittedvalues, series)
MSE2=mean_squared_error(fit2.fittedvalues, series)
MSE3=mean_squared_error(fit3.fittedvalues, series)
MSE4=mean_squared_error(fit4.fittedvalues, series)

#=====================================
# Printing the paramters and errors for each scenario
#=====================================
results=pd.DataFrame(index=[r"alpha", r"beta", r"gamma", r"l0", "b0", "MSE"])
params = ['smoothing_level', 'smoothing_slope', 'smoothing_seasonal', 'initial_level', 'initial_slope']
results["HW model 1"] = [fit1.params[p] for p in params] + [MSE1]
results["HW model 2"] = [fit2.params[p] for p in params] + [MSE2]
results["HW model 3"] = [fit3.params[p] for p in params] + [MSE3]
results["HW model 4"] = [fit4.params[p] for p in params] + [MSE4]
print(results)

#=====================================
# Evaluating and plotting the residual series for each scenario
#=====================================
residuals1= fit1.fittedvalues - series
residuals2= fit2.fittedvalues - series
residuals3= fit3.fittedvalues - series
residuals4= fit4.fittedvalues - series
residuals1.rename('residual plot for model 1').plot(color='red', legend=True)
residuals2.rename('residual plot for model 2').plot(color='blue', legend=True)
residuals3.rename('residual plot for model 3').plot(color='green', legend=True)
residuals4.rename('residual plot for model 4').plot(color='yellow', legend=True)
pyplot.title('Residual plots for models 1-4')
pyplot.show()

#=====================================
# ACF plots of the residual series for each scenario
#=====================================
from statsmodels.graphics.tsaplots import plot_acf
plot_acf(residuals1, title='Residual ACF for model 1', lags=50)
plot_acf(residuals2, title='Residual ACF for model 2', lags=50)
plot_acf(residuals3, title='Residual ACF for model 3', lags=50)
plot_acf(residuals4, title='Residual ACF for model 4', lags=50)
pyplot.show()
\end{lstlisting}

\section{Codes on ARIMA methods}
\begin{lstlisting}[language=Python, caption={WhiteNoiseModel}, label={WhiteNoiseModel}]
from random import gauss
from random import seed
from pandas import Series
from statsmodels.graphics.tsaplots import plot_acf, plot_pacf
from matplotlib import pyplot

# seed random number generator
seed(1)
# create white noise series
series = [gauss(0.0, 1.0) for i in range(1000)]

# Once created, we can wrap the list in a Pandas Series for convenience.
series = Series(series)

# summary statistics of the artificially generated series
print('Statistics of the artificially generated series:')
print(series.describe())

# line plot of the artificially generated series
series.plot(title='Time plot of a white noise model')
pyplot.show()

# histogram plot of the artificially generated series
series.hist()

# ACF plot of an artificially generated white noise time series
plot_acf(series, title='ACF of a white noise model', lags=50)

# PACF plot of an artificially generated white noise time series
plot_pacf(series, title='PACF of a white noise model', lags=50)
pyplot.show()
\end{lstlisting}

\begin{lstlisting}[language=Python, caption={AcfPacfPlotNonStationarity.py}, label={AcfPacfPlotDowJones.py}]
from pandas import read_excel
from matplotlib import pyplot
from statsmodels.graphics.tsaplots import plot_acf, plot_pacf
series = read_excel('DowJones.xls', sheet_name='Data2', header=0, index_col=0, parse_dates=True, squeeze=True)

# ACF plot on 50 time lags
plot_acf(series, title='ACF of Dow Jones time series', lags=50)

# PACF plot on 50 time lags
plot_pacf(series, title='PACF of Dow Jones time series', lags=50)
pyplot.show()
\end{lstlisting}

\begin{lstlisting}[language=Python, caption={ADFtest.py}, label={ADFtest.py}]
from pandas import read_excel
import matplotlib.pyplot as plt
series = read_excel('BuildingMaterials.xls', sheetname='Data', header=0,
              index_col=0, parse_dates=True, squeeze=True)
series.plot(color='red')
plt.xlabel('Dates')
plt.ylabel('Production values')
plt.title('Building materials production from 1986 to 2008')
plt.show()
#------------------------
from pandas import read_excel
from matplotlib import pyplot
from statsmodels.graphics.tsaplots import plot_acf, plot_pacf
series = read_excel('BuildingMaterials.xls', sheet_name='Data', usecols = [1],
                             header=0, squeeze=True, dtype=float)
# ACF plot on 50 time lags
plot_acf(series, title='ACF of building materials time series', lags=50)

# PACF plot on 50 time lags
plot_pacf(series, title='PACF of building materials time series', lags=50)
pyplot.show()
#--------------------------------------------------
# calculate stationarity test of time series data
#from pandas import read_csv
from statsmodels.tsa.stattools import adfuller
#series = read_csv('airline-passengers.csv', header=0, index_col=0, parse_dates=True,squeeze=True)
X = series.values
result = adfuller(X)
print('ADF Statistic: %f' % result[0])
print('p-value: %f' % result[1])
print('Critical Values:')
for key, value in result[4].items():
    print('\t%s: %.3f' % (key, value))
\end{lstlisting}

\begin{lstlisting}[language=Python, caption={AcfPacfAR1model.py}, label={AcfPacfAR1model.py}]
from pandas import read_excel
from statsmodels.graphics.tsaplots import plot_acf, plot_pacf
from matplotlib import pyplot

series = read_excel('DataAR1model.xls', sheet_name='ARdata', header=0, index_col=0, squeeze=True)

# Time, ACF, and PACF plots for original data
pyplot.plot(series)
pyplot.title('Time plot AR1 data')
plot_acf(series, title='ACF plot of AR1 data', lags=20)
plot_pacf(series, title='PACF plot of AR1 data', lags=20)
pyplot.show()
\end{lstlisting}

\begin{lstlisting}[language=Python, caption={AcfPacfMA1model.py}, label={AcfPacfMA1model.py}]
from pandas import read_excel
from statsmodels.graphics.tsaplots import plot_acf, plot_pacf
from matplotlib import pyplot

series = read_excel('DataMA1model.xls', sheet_name='MAdata', usecols = [1],
                    header=0, squeeze=True)

# Time, ACF, and PACF plots for original data
pyplot.plot(series)
pyplot.title('Time plot MA1 data')
plot_acf(series, title='ACF plot of MA1 data', lags=20)
plot_pacf(series, title='PACF plot of MA1 data', lags=20)
pyplot.show()
\end{lstlisting}

\begin{lstlisting}[language=Python, caption={AutomaticOrderARIMA.py}, label={AutomaticOrderARIMA.py}]
#===================================================
#Code for identifying the parameters with smallest AIC
#===================================================
from pandas import read_excel
from statsmodels.tsa.arima_model import ARIMA
import warnings
import itertools

series = read_excel('BuildingMaterials.xls', sheet_name='Data', header=0, index_col=0, parse_dates=True, squeeze=True)

#define the p, d and q parameters to take any value between 0 and 2
p = d = q = range(0, 3)

#generate all different combinations of p, q and q triplets
pdq = list(itertools.product(p, d, q))

#indentification of best model from different combinations of pdq
warnings.filterwarnings("ignore") # specify to ignore warning messages
best_score, best_param = float("inf"), None
for param in pdq:
        try:
            mod = ARIMA(series, param)
            results = mod.fit()
            if results.aic < best_score:
                best_score, best_param = results.aic, param
            print('ARIMA{} - AIC:{}'.format(param, results.aic))
        except:
            continue
print('The best model is ARIMA{} - AIC:{}'.format(best_param, best_score))
\end{lstlisting}

\begin{lstlisting}[language=Python, caption={ARIMA.py}, label={ARIMA.py}]
from pandas import read_excel
import matplotlib.pyplot as plt
from statsmodels.tsa.arima_model import ARIMA

series = read_excel('BuildingMaterials.xls', sheet_name='Data', header=0,
                index_col=0, parse_dates=True, squeeze=True)

#ARIMA(1,1,2) model
model = ARIMA(series, order=(1,1,2))

# generates ARIMA Model Results table
model_fit = model.fit(disp=0)
fitted = model.fit(disp=1)
print(model_fit.summary())

# Actual vs Fitted
model_fit.plot_predict(dynamic=False)
plt.show()

# Forecast
fc, se, conf = fitted.forecast(15, alpha=0.005)  # 95% conf

pred_ci = pred_uc.conf_int()
# plotting forecasts ahead
ax = series.plot(label='Original data')
pred_uc.predicted_mean.plot(ax=ax, label='Forecast values', title='Forecast plot with confidence interval')
ax.fill_between(pred_ci.index,
                pred_ci.iloc[:, 0],
                pred_ci.iloc[:, 1], color='k', alpha=.25)
plt.legend()
plt.show()
#====================================================================

#====================================================================
# MSE evaluation
y_forecasted = pred.predicted_mean
y_truth = series['2000-01-01':]
# Compute the mean square error
MSE = ((y_forecasted - y_truth) ** 2).mean()
print('MSE of the forecasts is {}'.format(round(MSE, 2)))
\end{lstlisting}

\begin{lstlisting}[language=Python, caption={FirstPlusSeasonalDifference.py}, label={FirstSeasonalDifference.py}]
#removing seasonality from time series using seasonal differencing
from pandas import read_excel
from statsmodels.graphics.tsaplots import plot_acf, plot_pacf
from matplotlib import pyplot

series = read_excel('Electricity.xls', sheet_name='ELEC', header=0, index_col=0, parse_dates=True, squeeze=True)

#time, ACF, and PACF plots for original data
pyplot.plot(series)
pyplot.title('Time plot original data')
plot_acf(series, title='ACF plot of original data', lags=50)
plot_pacf(series, title='PACF plot of original data', lags=50)
pyplot.show()

#seaonal difference
X = series.values
SeasDiff = list()
for i in range(12, len(X)):
	value = X[i] - X[i - 12]
	SeasDiff.append(value)

#time, ACF, and PACF plots for the seasonally differenced series
pyplot.plot(SeasDiff)
pyplot.title('Time plot seasonally differenced series')
plot_acf(SeasDiff, title='ACF plot of seasonally differenced series', lags=50)
plot_pacf(SeasDiff, title='PACF plot of seasonally differenced series', lags=50)
pyplot.show()

#seasonal + first difference
Y = SeasDiff
SeasFirstDiff = list()
for i in range(1, len(Y)):
	value = Y[i] - Y[i - 1]
	SeasFirstDiff.append(value)
pyplot.plot(SeasFirstDiff)
pyplot.title('Time plot seasonally + first differenced series')
plot_acf(SeasFirstDiff, title='ACF plot of seasonally + first differenced series', lags=50)
plot_pacf(SeasFirstDiff, title='PACF plot of seasonally + first differenced series', lags=50)
pyplot.show()
\end{lstlisting}

\begin{lstlisting}[language=Python, caption={AutomaticOrderSARIMAX.py}, label={AutomaticOrderSARIMAX.py}]
#===================================================
#Code for identifying the parameters with smallest AIC
#===================================================
from pandas import read_excel
import matplotlib.pyplot as plt
import statsmodels.api as sm
import warnings
import itertools
plt.style.use('fivethirtyeight')

series = read_excel('BuildingMaterials.xls', sheet_name='Data', header=0, index_col=0, parse_dates=True, squeeze=True)

#Define the p, d and q parameters to take any value between 0 and 1
p = d = q = range(0, 2)

# Generate all different combinations of p, q and q triplets
pdq = list(itertools.product(p, d, q))

# Generate all different combinations of seasonal p, q and q triplets (i.e., P, D, Q)
seasonal_pdq = [(x[0], x[1], x[2], 12) for x in list(itertools.product(p, d, q))]

# Indentification of best model from different combinations of pdq and seasonal_pdq
warnings.filterwarnings("ignore") # specify to ignore warning messages
best_score, best_param, best_paramSeasonal = float("inf"), None, None
for param in pdq:
    for param_seasonal in seasonal_pdq:
        try:
            mod = sm.tsa.statespace.SARIMAX(series, order=param, seasonal_order=param_seasonal, enforce_invertibility=False)
            results = mod.fit()
            if results.aic < best_score:
                best_score, best_param, best_paramSeasonal = results.aic, param, param_seasonal
            print('ARIMA{}x{} - AIC:{}'.format(param, param_seasonal, results.aic))
        except:
            continue
print('The best model is ARIMA{}x{} - AIC:{}'.format(best_param, best_paramSeasonal, best_score))
\end{lstlisting}

\begin{lstlisting}[language=Python, caption={SARIMAX.py}, label={SARIMAX.py}]
from pandas import read_excel
import pandas as pd
import matplotlib.pyplot as plt
import statsmodels.api as sm
plt.style.use('fivethirtyeight')

#==================================================================
#Loading the data set
df = read_excel('BuildingMaterials.xls', sheet_name='Data', header=0,
                index_col=0, parse_dates=True, squeeze=True)
#==================================================================

#==================================================================
#Fitting the ARIMA model and printing related statistics
# ARIMA(0, 1, 1)(0,1,1)12 in this case;
#this one is based on MA1 model template
mod = sm.tsa.statespace.SARIMAX(df, order=(1,1,1), seasonal_order=(0,1,1,12))
results = mod.fit(disp=False)
print(results.summary())
#==================================================================

#GRAPH BLOCK1======================================================
#Printing the graphical statistics of model (correlogram = ACF plot)
results.plot_diagnostics(figsize=(15, 12))
plt.show()
#==================================================================

#GRAPH BLOCK2======================================================
# printing the part of forecasts fitted to original data (for accuracy evaluation)
# the start date has to be provided; hence should be a time within the original time series;
# in this case, it is to start on 01 January 2000
pred = results.get_prediction(start=pd.to_datetime('2000-01-01'), dynamic=False)
pred_ci = pred.conf_int()

# printing one-step ahead forecasts together with the original data set;
# hence, the starting point (year) of the data set is required
# in order to build the plot of original series
ax = df['1986':].plot(label='Original data')
pred.predicted_mean.plot(ax=ax, label='One-step ahead Forecast', alpha=.7)
ax.fill_between(pred_ci.index,
                pred_ci.iloc[:, 0],
                pred_ci.iloc[:, 1], color='k', alpha=.2)
plt.legend()
plt.show()
#===================================================================

#GRAPH BLOCK3=======================================================
# Get forecast 20 steps ahead in future
pred_uc = results.get_forecast(steps=20)
# Get confidence intervals of forecasts
pred_ci = pred_uc.conf_int()
# plotting forecasts ahead
ax = df.plot(label='Original data')
pred_uc.predicted_mean.plot(ax=ax, label='Forecast values', title='Forecast plot with confidence interval')
ax.fill_between(pred_ci.index,
                pred_ci.iloc[:, 0],
                pred_ci.iloc[:, 1], color='k', alpha=.25)
plt.legend()
plt.show()
#====================================================================

#====================================================================
# MSE evaluation
y_forecasted = pred.predicted_mean
y_truth = df['2000-01-01':]
# Compute the mean square error
MSE = ((y_forecasted - y_truth) ** 2).mean()
print('MSE of the forecasts is {}'.format(round(MSE, 2)))
#====================================================================
\end{lstlisting}

\section{Codes on the regression analysis and application to forecasting}

%\begin{lstlisting}[language=Python, caption={RegressionTimePlot.py}, label={RegressionTimePlot.py}]
%from pandas import read_excel
%import matplotlib.pyplot as plt
%series = read_excel('Bank.xls', sheet_name='Data2', header=0,
%                     squeeze=True, dtype=float)
%
%#reading the basic variables
%DEOM = series.DEOM
%AAA = series.AAA
%Tto4 = series.Tto4
%D3to4 = series.D3to4
%
%DEOM.plot()
%plt.xlabel('time')
%plt.ylabel('Difference end of month balance')
%plt.title('DEOM')
%plt.show()
%
%AAA.plot()
%plt.xlabel('time')
%plt.ylabel('Composite AAA Bond rates')
%plt.title('AAA')
%plt.show()
%
%Tto4.plot()
%plt.xlabel('time')
%plt.ylabel('US Govt 3-4 year Bond rates')
%plt.title('3to4')
%plt.show()
%
%D3to4.plot()
%plt.xlabel('time')
%plt.ylabel('Difference US Govt 3-4 year Bond rates')
%plt.title('D3to4')
%plt.show()
%\end{lstlisting}

\begin{lstlisting}[language=Python, caption={RegressionAnalysis.py}, label={Regression.py}]
from pandas import read_excel
from statsmodels.formula.api import ols
series = read_excel('Bank.xls', sheet_name='Data2', header=0,
                     squeeze=True, dtype=float)

#reading the basic variables
DEOM = series.DEOM
AAA = series.AAA
Tto4 = series.Tto4
D3to4 = series.D3to4

#Regression model(s)
formula = 'DEOM ~ AAA + Tto4 + D3to4'

#Ordinary Least Squares (OLS)
results = ols(formula, data=series).fit()
print(results.summary())

# Here the main table is the first one,
# where the main statistics are the R-squared (line 1)
# and the P-value; i.e., Prob (F-statistic)
\end{lstlisting}

\begin{lstlisting}[language=Python, caption={RegressionForecast.py}, label={RegressionForecast.py}]
from pandas import read_excel
from statsmodels.tsa.api import Holt
#from matplotlib import pyplot
import matplotlib.pyplot as plt
import numpy as np
from statsmodels.formula.api import ols
series = read_excel('Bank.xls', sheet_name='Data2', header=0,
                     squeeze=True, dtype=float)

# Reading the basic variables
DEOM = series.DEOM
AAA = series.AAA
Tto4 = series.Tto4
D3to4 = series.D3to4

# Forecasting for AAA using Holt's linear method
fit1 = Holt(AAA).fit(optimized=True)
fcast1 = fit1.forecast(6).rename("Additive 2 damped trend")
fit1.fittedvalues.plot(color='red')
fcast1.plot(color='red', legend=True)
AAA.plot(color='black', legend=True)
plt.title('Forecast of AAA with Holt linear method')
plt.show()

# Forecasting for Tto4 using Holt's linear method
fit2 = Holt(Tto4).fit(optimized=True)
fcast2 = fit2.forecast(6).rename("Additive 2 damped trend")
fit2.fittedvalues.plot(color='red')
fcast2.plot(color='red', legend=True)
Tto4.plot(color='black', legend=True)
plt.title('Forecast of 3to4 with Holt linear method')
plt.show()

# Forecasting for D3to4 using Holt's linear method
fit3 = Holt(D3to4).fit(smoothing_level=0.8, smoothing_slope=0.2, optimized=False)
#fit3 = Holt(D3to4).fit(optimized=True)
fcast3 = fit3.forecast(6).rename("Additive 2 damped trend")
fit3.fittedvalues.plot(color='red')
fcast3.plot(color='red', legend=True)
D3to4.plot(color='black', legend=True)
plt.title('Forecast of D3to4 with Holt linear method')
plt.show()

# Building the regression based forecast for main variable, DEOM
# Regression model(s)
formula = 'DEOM ~ AAA + Tto4 + D3to4'

# ols generate statistics and the parameters b0, b1, etc., of the model
results = ols(formula, data=series).fit()
results.summary()
b0 = results.params.Intercept
b1 = results.params.AAA
b2 = results.params.Tto4
b3 = results.params.D3to4

# putting the fitted values of the forecasts of AAA, Tto4, and D3to4 in arrays
a1 = np.array(fit1.fittedvalues)
a2 = np.array(fit2.fittedvalues)
a3 = np.array(fit3.fittedvalues)

# Building the fitted part of the forecast of DEOM
F=a1
for i in range(53):
    F[i] = b0 + a1[i]*b1 + a2[i]*b2 + a3[i]*b3

# putting the values of the forecasts of AAA, Tto4, and D3to4 in arrays
v1=np.array(fcast1)
v2=np.array(fcast2)
v3=np.array(fcast3)

# Building the 6 values of the forecast ahead
E=v1
for i in range(6):
    E[i] = b0 + v1[i]*b1 + v2[i]*b2 + v3[i]*b3


# Joining the fitted values of the forecast and the points ahead
K=np.append(F, E)

# Reading the original DEOM time series for all the 59 periods
DEOMfull0 = read_excel('Bank.xls', sheet_name='Data4', header=0,
                     squeeze=True, dtype=float)

###########################
# Evaluating the MSE to generate the confidence interval
DEOMfull = DEOMfull0.DEOMfull
values=DEOMfull[0:53]
Error = values - F
MSE=sum(Error**2)*1.0/len(F)

## Lower and upper bounds of forecasts for z=1.282; see equation (2.2) in Chap 2.
#LowerE = E - 1.282*MSE
#UpperE = E + 1.282*MSE

LowerE = DEOMfull0.LowerE
UpperE = DEOMfull0.UpperE

print(LowerE)
###############################

# Plotting the graphs of K and DEOMfull with legends
from matplotlib.legend_handler import HandlerLine2D
line1, = plt.plot(K, color='red', label='Forecast values')
line2, = plt.plot(DEOMfull, color='black', label='Original data')
line3, = plt.plot(LowerE, color='blue', label='Lower forecast')
line4, = plt.plot(UpperE, color='orange', label='Upper forecast')
plt.legend(handler_map={line1: HandlerLine2D(numpoints=4)})
plt.title('DEOM regression forecast with confidence interval')
plt.show()

# Proceeding as as in other demos, forecasts
# can be generated for other scenarios; i.e., with different combinations of variables
\end{lstlisting}

\end{document}